 \documentclass[draft]{article}

\usepackage{amsmath,amsfonts,amsthm,amssymb,amscd,cancel,color,mathtools}
\usepackage{enumitem}
\usepackage{verbatim}
\usepackage{ulem,color}

\renewcommand{\Re}{\mathop{\rm Re}\nolimits}

\def\uno{{\kern+.3em {\rm 1} \kern -.22em {\rm l}}}

\theoremstyle{plain} \newtheorem{theorem}{Theorem}[section]
\newtheorem{lemma}[theorem]{Lemma}
\newtheorem{proposition}[theorem]{Proposition}
 \theoremstyle{definition}
\newtheorem{definition}[theorem]{Definition} \theoremstyle{remark}
\newtheorem{remark}[theorem]{Remark}

\newcommand{\R}{{\mathbb R}}
\newcommand{\U}{{\mathcal U}}

\newcommand{\Z}{{\mathbb Z}}

\newcommand{\N}{{\mathbb N}}

\def\im{{\rm i}}

\newcommand{\C}{\mathbb{C}}

\def\({\left(}
\def\){\right)}
\def\<{\left\langle}
\def\>{\right\rangle}
\newcommand{\rad}{{\mathrm{rad}}}
\def\[{\left [}
\def\]{\right ]}

\numberwithin{equation}{section}

\begin{document}

\title{A survey on asymptotic stability of ground states of nonlinear Schr\"odinger equations II }

\author {Scipio Cuccagna, Masaya Maeda }

\maketitle

\begin{abstract} We give short survey on the question of asymptotic stability of ground states of nonlinear Schr\"odinger equations,
 focusing primarily on the so called nonlinear  Fermi Golden Rule.

\end{abstract}

\section{Introduction}
\label{sec:introduction}

In  2004 one of us authored a survey \cite{cuccagna2004} on the asymptotic stability of ground states of the  nonlinear
Schr\"odinger equation (NLS). Since then there has been
considerable progress on this topic, so that it is worthwhile to write a review with some updates.

For $d\geq 1$, we consider the    NLS
\begin{equation}\label{NLSsingle}
 \im \partial_t u =     -\Delta u   +\beta  (|u|^2) u   , \quad \left .u\right | _{t=0}=u_0 \in H^1 (\R ^d, \C ),
\end{equation}
where $\beta \in C^\infty(\R,\R) $ satisfies, for $d^*=\infty$ for $d=1,2$ and $d^*=\frac{d+2}{d-2}$  for $d\ge 3$,
\begin{align}\label{eq:growth_condition}
\text{ $|   \partial ^{ n}_t  (\beta  (t^2) t)   | \le C_n t^{p-n} $ for $t\ge 1$, $n\ge 0$
	and for a  $p<d^* $.}
\end{align}
This guarantees that the Cauchy problem  \eqref{NLSsingle}  is locally well posed, see Cazenave  \cite{caz}.

We are concerned with a spatially localized solution called soliton.
In particular, we   assume there exists an open interval $\mathcal{O}\subset (0,\infty)$ such that
\begin{equation}
\label{eq:B}
\Delta u -\omega u-\beta(|u|^2)u=0\quad\text{for $x\in \R^d$},
\end{equation}
admits a $C^\infty$-family of ground states $\mathcal{O}\ni \omega \mapsto \phi _ {\omega } \in H ^{1}_{\rad}(\R  ^d)$  with $ \phi _ {\omega }(x) >0$  everywhere.
In fact, under these hypotheses, we have $\phi _ {\omega } \in C ^{\infty}(\R  ^d)$ and
\begin{equation}\label{eq:dec sol}
| \partial ^{\alpha}_x \phi _ {\omega } (x)|\le C _{\alpha , \omega} (1+ |x|) ^{-\frac{d-1}{2}}    e ^{-\sqrt{\omega}|x|} \text{ for all multiindexes }\alpha .
\end{equation}
Then $ e^{\im \(\frac 1 2 v\cdot x
	-\frac 1 4 |v|^2 t +  t\omega +  \vartheta \) }{ \phi} _\omega
(x-vt- D) $, for any choice of $(\omega , \vartheta , v ,D)\in \mathcal{O}\times \R \times \R ^d\times \R ^d$, are   solitonic solutions  of the NLS.
An important question is whether these ground state solutions are stable. A first notion of stability is the following.

\begin{definition}[Orbital stability]\label{def:orbst}
	A ground  state $   \phi _\omega
	$ of
	\eqref{NLS}  is {\it orbitally stable} if
	\begin{equation*}\text{$\forall$ $\epsilon>0$,  $\exists$ $\delta>0$ s.t.\    $\| {\phi }_\omega - {u}_0\|_{H^1}<\delta$}\Rightarrow
	\sup_{t>0}\inf_{(\vartheta ,D)\in\R \times \R^d}\|e^{\im   \vartheta } {\phi} _\omega (\cdot + D) - {u}(t)\|_{H^1}<\epsilon ,
	\end{equation*}
	where $ {u}$  is the solution of the NLS  with $ {u}(0)= {u}_0$.
\end{definition}
The literature on this is large,  see the survey papers \cite{BGR15Survey,Stuart08} and the  references therein.

\begin{theorem}\label{thm:class orb stab}
	If  for  $\omega \in\mathcal{O}$ both the two  conditions (H1), (H2) listed below  are satisfied, then the corresponding ground state is orbitally stable:
	\begin{itemize} \item[(H1)]   $\ker L _{+,\omega} \cap H^1 _{\rad}(\R ^d)=\{ 0\}$   for the    operator       $L_{+,\omega}:= -\Delta +\omega + \beta  (\phi _\omega ^2)+2\beta ' (\phi _\omega ^2)\phi _\omega ^2 $;
		\item[(H2)]   we have the  Vakhitov--Kolokolov condition  $\displaystyle  q'(\omega )  >0$, where $ q (\omega ) :=Q(\phi _\omega ) $.
	\end{itemize}
\end{theorem}

The  study of equilibria and of  solitons of NLS's or of  more complex models  and their orbital stability is not the topic of this paper.
The notion of orbital stability applies also to other functions. For example, if  $B(t)\ge 0$ for $t\ge 0$,
then energy and mass conservation and Gagliardo--Nirenberg inequalities imply in an elementary fashion the orbital stability of the 0 solution.
Less elementary is the following fact.   If 0 is stable and
\begin{equation}\label{eq:beh in 0}
\text{$d\ge 2$ and  for $d=1$  furthermore $\beta '(0) =0$,}
\end{equation}
then there is an $\epsilon _0>0$
s.t.
\begin{align}
\| {u} _0\| _{H^1(\R ^d)}<\epsilon _0 \Longrightarrow \exists \  {u}_+\in H^1(\R ^d) \text{ s.t.\ } \| {u}(t )-e^{\im  t \Delta }  {u}_+    \|_{H^1} \xrightarrow{t\to +\infty}0. \label{eq:as_stab_0}
\end{align}
The theme of the present paper is an analogue of \eqref{eq:as_stab_0} in the case of solitons of the NLS. Specifically,   we will give an outline of some of
the most basic ideas behind the   following analogous rough statement, which we call the asymptotic stability of solitons.

\begin{theorem}\label{thm:main rough}
	Let $d\geq 3$.
	Let $\omega _1 \in\mathcal{O}$ satisfy the two conditions listed in Theorem \ref{thm:class orb stab}. Then, under
	further hypotheses, which the authors of this review believe to  hold generically,  there exists $\epsilon _1>0$ s.t.\ for any $u_0\in B _{H^1   } (\phi _{\omega _1}   , \epsilon _1):=\{v\in H^1\ |\ \|v-\phi_{\omega _1} \|< \epsilon_1 \} $
	there exist $\omega _+\in \mathcal{O}$, $v_+\in \R ^d$
	and  $(\vartheta , D) \in C^0([0,+\infty ), \R \times \R^d)
	$ s.t.\ the solution of the NLS with  $u(0)=u_0$ satisfies
	\begin{equation}\label{eq:scattering}
	\| {u}(t )-e^{\im   \vartheta(t)+  \frac  \im  2{v_+ \cdot x} }   {\phi}_{\omega
		_+} (\cdot -D(t))- e^{\im   t\Delta }h _+    \|_{H^1(\R ^d) } \xrightarrow{t\to +\infty}0.
	\end{equation}
\end{theorem}

\begin{remark} \label{rem:dim1} For dimensions 1 and 2 the same theorem is known to be true  only  under conditions that break  translation, as when   $u_0$ is an even function
	or there is an additional translation breaking inhomogeneity in \eqref{NLS}, like a linear potential. The proof in the case with translation
	is an open problem.
\end{remark}

\begin{remark} \label{rem:dim1integr}  In dimension 1, well known is the case when $\beta  (|u|^2) u=-|u|^2  u $, where it is possible to apply methods from the theory
of integrable systems \cite{borghese,saalmann}, which require $u_0$ s.t. $\< x \> ^{s}u_0\in L^2(\R )$  for $s>1/2$, see \cite{cupe2014}.
\end{remark}

\begin{remark} \label{rem:generic1} The additional hypotheses required are (H3) (see Theorem \ref{thm:BP1}), (H4)--(H7) in Sect. \ref{sec:eigenvalues} and (H8) under
	\eqref{eq:FGR000}.  The most delicate condition in  (H8)    requires that the terms in \eqref{eq:FGR000} be non zero. This happens when the  Fourier transform
	of certain functions has nonzero restriction on   certain spheres  of phase space.   When $\beta $ is real analytic, then the dependence of the coefficients on  $\omega$ is analytic.
	
	\noindent  The generic condition has not been proved rigorously, except in very special situations, see \cite{Komech2012,Buslaev2008,adami13JMP}.  Even the question of checking  numerically the generic condition
	seems to have attracted very little interest.
\end{remark}

%
%
	Theorem \ref{thm:main rough} has a long history.
	The theory was initiated by Soffer and Weinstein \cite{SW1,SW2} for small solitons bifurcating from linear potential, see also \cite{PW},  followed by important paper by Buslaev and Perelman \cite{BP1} which proved the asymptotic stability for the case $d=1$.
	Both of \cite{SW1,SW2} and \cite{BP1} considers the case where the linearized operator have no non-zero eigenvalues.
	The basic idea of these works is to divide the solution into a soliton part and a remainder part by modulation argument and then prove the decay of the remainder part by the dispersive properties of the linearized operator. The remainder is small and satisfies a complicated equation that looks like a NLS. The linear part of the equation of the remainder, has continuous spectrum and eigenvalues.

Very early  the literature provided a theory of the dispersive properties of the continuous part of the linearized operator. In dimension $d=1$ this is in \cite{BP1}, which can be supplemented with Krieger and Schlag \cite{KS}, see also \cite{zhousigal2}
the case   $d\ge 3 $  is in  \cite{Cu1}, which has to be supplemented with  \cite{CPV}, and $d=2$ in \cite{cuccagnatarullidim2}. More effective use   of dispersion, of Strichartz estimates and especially of the endpoint Strichartz estimate in $d\ge 3$, see Keel and Tao  \cite{Kl-Tao}, is  in Gustafson, Nakanishi and Tsai \cite{GNT}.
	Smoothing estimates as a surrogate of the   endpoint Strichartz estimate when  $d=1,2$
	are in Mizumachi
	\cite{M1,M2}. A substantial simplification  of  Mizumachi's smoothing estimates is in   \cite{cuccagnatarulliDNLS}. Other early contributions are \cite{kirr1,kirr2}.     Obviously, dispersion is a hard problem in the presence
of strong nonlinearities, where one cannot hope to prove dispersive properties of the remainder  just by Strichartz estimates, and here the literature is not as rich. Remarkable nonetheless are   \cite{KMM2,KMM3} as well as the very recent  \cite{KMM4}.

While, to some extent, linear dispersion of the continuous mode was understood quite early, it took some time to understand   how to treat the nonzero eigenvalues of the linearized operator.  The starting point  seems to be  Sigal \cite{Sigal93CMP} which, for a different problem, showed
the existence of a nonlinear damping mechanism  by which the discrete modes lose energy which, by nonlinear interaction, spills in the continuous part of the equation and then scatters by essentially linear mechanisms.
Sigal called this damping mechanism "nonlinear Fermi Golden Rule" (FGR).
The first successful implementation of this idea in our context was obtained by Buslaev and Perelman   \cite{BP2}. For almost 15 years there was no major improvement on this part of
the proof in \cite{BP2}. Here we recall  that \cite{BP2} treats the case where there is just a single     $\mathbf{e}(\omega) \in (0,\omega )$
of multeplicity 1 of the linearization  operator $\mathcal{H}_{\omega}$,  with $2\mathbf{e}(\omega)> \omega$. Later   Soffer and Weinstein \cite{SW3} developed a similar idea
in the context of the NLKG equation.   See also \cite{PKA98}.
 Various papers where written in the early '00 \cite{BS,TY1,TY2,TY3,T,Cu2,SW4,zhouweinstein1,zhouweinstein2} articulating the idea. A novelty was  in  Gang Zhou and Sigal \cite{zhousigal} , with still
   just one eigenvalue but with $2$ replaced by $N+1$ for  $N\in \N$, see also   \cite{cuccagnamizumachi}. However, all these rather restrictive conditions on the spectrum of the linearized operator $\mathcal{H}_{\omega}$, where finally lifted   only with  \cite{bambusicuccagna,Cuccagna11CMP} around 2010.  These papers introduced a more natural framework for a problem that, approached from a different viewpoint,
   could look  impossibly complex, as can be seen, for example, by tracing the argument in \cite{Gz}. It should be remarked, that quite independently from the theory we are discussing here,
   Perelman  \cite{perelman01} and Merle and Raphael \cite{MR4} exploited a form of FGR in their masterly   analysis of the    $\| \nabla u(t)\| _{L^2(\R ^d)}\sim \frac{\sqrt{\log |\log  t |}}{t}$  blow up
in the NLS with $\beta  (|u|^2) u=-|u|^{\frac{4}{d}}u$.  The connections between the two theories have not been explored yet, although \cite{CM20} exploits ideas originating from the work of Merle and Raphael to simplify considerably the proof of the result in \cite{CM15APDE}.

The paper \cite{Cuccagna11CMP}
considers equations without translation. Translation was later and  independently  introduced in  \cite{Cuccagnatrans} and \cite{Bambusi13CMPas}.
However, there are aspects of the proof, which is rather long and with many detains, that have been finalized in later papers, such as \cite{CM1}. See also \cite{BM16CMP}
for some more on \cite{Bambusi13CMPas}.

	\noindent In this paper  we will just focus on the FGR.
	The generic conditions in  {Theorem} \ref{thm:main rough} pertain to the FGR. As we mentioned, there is very little numerical work on them.
	The fact that a certain quadratic form is non--negative, is explained  later. Strict positivity is unproven, theoretically as well as numerically.  Numerical simulations are certainly not simplified by the fact
	that the crucial quadratic form is obtained after a rather complex sequence of coordinate changes. The coordinate changes are not discussed in this paper.
	

The aim of this survey is to give some basic intuition of the main ideas of the proof of {Theorem} \ref{thm:main rough}, skipping completely on the most technical parts of the proof.

\section{Theorem \ref{thm:main rough} in the absence of nonzero eigenvalues}
\label{sec:no eigenvalues}

We  embed $\C \hookrightarrow \C^2$
  using the natural identification
\begin{equation}\label{eq:identif}
  \C \ni u \mapsto \widetilde{u}:= \begin{pmatrix}
  u\\ \bar{u}
  \end{pmatrix}\in \widetilde{\C}:=\left \{ \begin{pmatrix}
  z\\ \bar{z}
  \end{pmatrix}\in \C^2 : z\in \C  \right \} \subset \C^2.
\end{equation}
Here  we set $\<U ,V\>_{\C^2}:=
2^{-1} (u_1  v_1+u_2   v_2 ) $ for  $U={}^t(u_1\   u_2)  $ and $V={}^t(v_1\   v_2) $ in $\C^2$.  By this definition, $\< \widetilde{u},\sigma _1\widetilde{v}\>_{\C^2}=\Re u\bar v$, and in particular $\<\widetilde{u},\sigma _1\widetilde{u}\>_{\C^2}=|u|^2$, where \begin{align}\label{17.3}  \sigma_1=\begin{pmatrix} 0 & 1 \\ 1 & 0 \end{pmatrix},\ \sigma_2=\begin{pmatrix} 0 & -\im \\ \im & 0 \end{pmatrix},\  \sigma_3=\begin{pmatrix} 1 & 0 \\ 0 & -1 \end{pmatrix}.\end{align}
Armed with this, we can equivalently write the NLS as
 \begin{align}\label{NLS}
\im \sigma_3 \partial_t \widetilde u =-\Delta \widetilde u + \beta ( \< \widetilde u ,\sigma _1\widetilde u \>_{\C^2}  )\widetilde u     , \quad \left .\widetilde{u}\right | _{t=0}=\widetilde{u}_0 \in H^1 (\R ^d, \widetilde{\C}  ).
\end{align}

  \noindent With the above definition of $\<\cdot  ,\cdot \>_{\C^2}$,  we define  \begin{equation}\label{eq:bilform}
   \<U,V\>:=\int_{\R^3}\<U(x),V(x)\>_{\C^2}\,dx \text{ for $U,V\in L^2(\R^3,\C^2)$  }
\end{equation}
(we emphasize, that here there is no complex conjugation).  In $L^2(\R ^d, \widetilde{\C} )$  we consider the symplectic form
$\Omega$, defined by
 \begin{equation}\label{eq:SymplecticForm}
  \Omega (X,Y)=\im \langle X, \sigma _3\sigma _1 Y \rangle  \text{  for all $X,Y\in L^2(\R ^d, \widetilde{\C} )$.}
\end{equation}
Given a function $F\in C^1( U,\R )$, with $U$ an open subset of $H^1(\R ^d, \widetilde{\C} )$, we denote by $dF(u)$ the Frech\' et derivative of
$F$, and  by $ \nabla F(\widetilde{u})$ its gradient, defined by  $dF(\widetilde{u}) =\< \sigma _1\nabla F(\widetilde{u}) , \cdot \>$.
The Hamiltonian vector--field  $X_F$   of $F$ associated to $ \Omega$  is defined by $ \Omega (X_F,\cdot )=dF$, that is $X_F =-\im \sigma _3 \nabla F$.

\noindent If we consider, for $B(0)=0$  the primitive $B'   =\beta  $, the energy
\begin{equation} \label{eq:energyfunctional}\begin{aligned}&
 E(\widetilde{u}):= \frac{1}{2}\<(-\Delta)\widetilde{u},\sigma _1\widetilde{u}\>+ \frac 1 2 \int_{\R^3}B(\< \widetilde u ,\sigma _1\widetilde u \>_{\C^2})\,dx,
  \end{aligned}
\end{equation}
     then $ \nabla E(\widetilde{u}) = -\Delta \widetilde u + \beta (\< \widetilde u ,\sigma _1\widetilde u \>_{\C^2})\widetilde u$
    and \eqref{NLS} can be interpreted as $\partial_t \widetilde u =X_E (\widetilde u )$.

    \noindent Notice that $E\in C^2 \( H^1(\R ^d, \widetilde{\C} ), \R \)$ with
    \begin{align}
        \nabla ^2 E(\widetilde{u})\widetilde{X}&:= \frac{d}{dt}  \nabla E(\widetilde{u}+t\widetilde{X})|_{t=0}= -\Delta \widetilde{X} +\beta (\< \widetilde u ,\sigma _1\widetilde u \>_{\C^2}) \widetilde{X}+ 2\beta '(\< \widetilde u ,\sigma _1\widetilde u \>_{\C^2}) \< \widetilde u ,\sigma _1\widetilde X \>_{\C^2}\widetilde{u} \label{eq:2nd_der}\\& = -\Delta \widetilde{X} +\beta (\< \widetilde u ,\sigma _1\widetilde u \>_{\C^2}) \widetilde{X} + \beta '(\< \widetilde u ,\sigma _1\widetilde u \>_{\C^2}) |u|^2\widetilde{X}
        + \beta '(\< \widetilde u ,\sigma _1\widetilde u \>_{\C^2}) \begin{pmatrix} 0 & u^2 \\ \overline{u}^2 & 0 \end{pmatrix}\widetilde{X}.\nonumber
    \end{align}

    \noindent We define also quadratic forms $P_j(u) :=  {2}^{-1}  \<  \Diamond _  j\widetilde{u},\sigma _1\widetilde{u} \>$ for $j=0,1,...,d$,
    which are invariant by  gauge and translation symmetries, with
  \begin{align}\label{eq:invariants_0}
 Q(u)  =P_0(u) \text{   for $\Diamond _0:=1$  the mass and  $ P_a(u)$ for $\Diamond _a :=-\sigma _3\im  \partial _{ a}$,     $a=1,...,d,$}
\end{align}
  the linear momenta.

Here we extend the hypotheses in {Theorem} \ref{thm:class orb stab}, and assume that  $   q'(\omega )  >0$ for all $\omega \in \mathcal{O}$, which can be assumed, if necessary, restricting $\mathcal{O}$.
  Under such assumption, the map $(\omega , v)\to  p=\Pi (e^{\sigma _3  \frac \im  2v\cdot x   }\phi _\omega     ) $ is a diffeomorphism  into an open
subset   $\mathcal{P}$ of $ \R ^d$.  This uses also $ \Pi _a (e^{ \sigma _3  \frac \im  2    v\cdot x   }u )= \Pi _a(u) +2^{-1} v_a Q (u
   ) $ for $a = 1,..., d$.
For $p=p(\omega, v)\in \mathcal{P}$ set $\Phi
_p=e^{ \sigma_3  \frac \im  2   v\cdot x   }\widetilde{ \phi} _\omega    $.
    The $\Phi _p   $ are
constrained critical points of $E$ with associated Lagrange
multipliers $\lambda  (p) \in \R^{d+1}$ so that
 \begin{equation}
	 \label{eq:eq sol} \nabla E( e^{ \im \sigma _3 \tau \cdot \Diamond}\Phi _p  )= \lambda   (p) \cdot \Diamond e^{ \im \sigma _3 \tau \cdot \Diamond}\Phi _p,
 \end{equation}
 where
 we have
 \begin{equation}
	 \label{eq:LagrMult} \lambda _0(p) =-\omega  (p) -4 ^{-1}{v^2 (p)}    \, ,    \quad  \lambda _a(p):=v_a (p)
	\, \text{ for $a=1,..., d$.}
 \end{equation}
We now introduce the linearization, for $(\omega, v)=(\omega (p), v(p))$,
\begin{align}\label{eq:linearizationL}&
 \mathcal {H}_p :=   \sigma _3 (\nabla ^2
E(\Phi _p  )- \lambda  (p) \cdot \Diamond )  =\sigma_3(-\Delta+ \omega    +4 ^{-1}{v^2 }  +\im v \cdot \nabla ) +V_{p}\\& \text{where }  V_{p}:=
   \sigma_3
\left[\beta (\phi ^2_{\omega (p) }) +\beta ^\prime (\phi ^2_{\omega (p)
})\phi ^2_{\omega (p) } \right] +\im  \sigma _2 \beta ^\prime (\phi ^2
_{\omega  (p)})\phi ^2 _{\omega (p) }    e^{-\sigma _3  \frac \im  2    v (p)\cdot x },  \nonumber
 \end{align}
 which can be computed from \eqref{eq:2nd_der}.
By an abuse of notation, we set
\begin{equation}\label{eq:linearization0}
   \mathcal {H}_\omega := \mathcal {H}_p  \text{ when $v(p)=0$ and $\omega (p)=\omega$.}
 \end{equation}
 It is easy that ${\mathcal H}_p =e^{  \sigma _3 \frac \im  2   v (p)\cdot x}
  {\mathcal H}_{\omega (p)} e^{ - \sigma _3 \frac \im   2   v (p)\cdot x}$, so that the spectrum of ${\mathcal H}_p$   depends only on $\omega(p)$.

  Hypothesis (H2)  of {Theorem} \ref{thm:class orb stab} guarantee that the map $p\to \lambda (p)$ is a local diffeomorphism and, in particular, it is invertible.
  In \cite{W2} it is shown that    Hypothesis (H1)  of {Theorem} \ref{thm:class orb stab}  implies  the following:  \begin{align}
&\ker {\mathcal H}_p  =\text{Span}\{  \sigma _3  \Diamond _j    \Phi _p:j=0,..., d   \} \text{ and} \label{eq:kernel1}\\
&N_g ( {\mathcal H}_p ) = \text{Span}\{ \sigma _3  \Diamond _j    \Phi _p, \partial _{ \lambda _j}     \Phi _p :j=0,..., d \} , \label{eq:kernel1-}
\end{align}
where   $N_g ( L ) :=\cup _{j=1}^\infty \ker (L^j)$.  Notice that the $\supseteq $  in  \eqref{eq:kernel1}   follows immediately differentiating in $\tau$ the identity
\eqref{eq:eq sol}  while the opposite inclusion is a much harder proposition,
 which rests  on  $\ker L_+ \cap H^1 _{\rad}(\R ^d)=\{ 0\}$. Setting $\tau=0$ in  \eqref{eq:eq sol} and differentiating in $\lambda _j$, we obtain
 the $\supseteq $  in  \eqref{eq:kernel1-}.  The $\subseteq $  in  \eqref{eq:kernel1-} follows from
 \eqref{eq:kernel1},
 the fact that the correspondence $p \longleftrightarrow \lambda$ is a diffeomorphism (this, in turn a consequence of  $   q'(\omega )  >0$ for all $\omega \in \mathcal{O}$),  from    Fredholm alternative and from
\begin{equation}\label{eq:orth-}
 \delta _{jk} =\partial _{p_k} p_j=   2^{-1} \partial _{p_k}  \< \Diamond _j \Phi _p, \Phi _p \> = \< \Diamond _j \Phi _p, \partial  _{p_k} \Phi _p \> .
\end{equation}
  We have the decomposition  \begin{align} 	
\label{eq:begspectdec2}& L^2(\R ^d, \C ^2) = N_g(\mathcal{H}_p)\oplus N_g^\perp
(\mathcal{H}_p^{\ast}) \  ,
   \\& N_g (\mathcal{H}_p^{\ast})  =\text{Span}\{   \Diamond _j    \Phi _p,
     \sigma _3 \partial _{ \lambda _j}     \Phi _p :j=0,..., d  \}   .
   \label{eq:begspectdec3}
\end{align}
Set  $P_{N_g }(p) =P_{N_g ( {\mathcal H}_p )}$ for the  projection on
$N_g(  {\mathcal H}_p )  $    and     $P(p):=1-P_{N_g }(p)$. Notice that
\begin{equation} \label{eq:projNg} \begin{aligned} & P_{N_g }(p)X=\sum _{j=0} ^{d}\(  \sigma _3 \Diamond
_j \Phi _p\  \langle \sigma _1X , \sigma _3  \partial _{p_j}\Phi _p\rangle +\partial
_{p_j}\Phi _p \  \langle \sigma _1X ,\Diamond _j \Phi _p \rangle \) . \end{aligned}
\end{equation}
Then we have the following Modulation Lemma, which originates with Soffer and Weinstein \cite{SW1}.

\begin{lemma}[Modulation]
  \label{lem:modulation}
     Fix  $p_1 \in \mathcal{P}$. Then there exists a neighborhood $\U   $   of $\Phi _{ p_1}$    in $ H^  {1}(\R ^d, \C ) $
    and   functions $p  \in C^\infty (\U     , \mathcal{P})$
	and $\tau \in C^\infty (\U    , \R ^{d+1  })$ s.t.\  $p(\Phi _{ p_1})=p_1 $   and $\tau
(\Phi _{ p_1}) =0$
  and s.t.\  $\forall u\in \U    $
\begin{equation}\label{eq:ansatz}
\begin{aligned}
  &
 u =   e^{\im \sigma _3  \tau \cdot \Diamond} (  \Phi _{p } +R)
 \text{  and $R\in N^{\perp}_g (\mathcal{L}_p ^*)$.}
\end{aligned}
\end{equation}

\end{lemma}
We write \eqref{NLS} as $\im \dot {\widetilde{u}}= \nabla E(\widetilde{u})$.
Using  \eqref{eq:ansatz},  $
 \nabla E(\Phi _p  )= \lambda   (p) \cdot \Diamond \Phi _p$ ,  the definition of ${\mathcal H}_p$ and for  $O(R^2)$ is non--linear in $R$, we obtain
\begin{align} &       - \( \dot \tau   -\lambda   (p) \)   \cdot   \sigma _3 \Diamond    \Phi _{p }     +\im\  \dot p \cdot   {\partial _p } \Phi _{p }-  \( \dot \tau   -\lambda   (p) \) \cdot  \sigma _3 \Diamond   R+  \im \dot R
 =     {\mathcal H}_p  R+O(R^2).\label{eq:eq ansatz}\end{align}
 Applying  $ P_{N_g }(p)$ to \eqref{eq:eq ansatz}, summing on repeated indexes,  we obtain the Modulation Equations
\begin{align}
  & \dot \tau _k   -\lambda _k+  \( \dot \tau  -\lambda    \)  \cdot  \<     \Diamond    R  ,   \partial _{p_k}\Phi _p   \>  -  \dot p  \cdot  \<   \im \sigma _1   R  ,\sigma _3 \partial _{p_k}\partial _{p}\Phi _p  \> = \<     O(R^2)  ,    \partial _{p_k}\Phi _p  \> \nonumber \\& \dot p _k   -  \( \dot \tau    -\lambda     \)  \cdot \<    \im   \sigma _3\sigma _1\Diamond     R  ,\Diamond _{k}\Phi _p   \>  - \dot p \cdot  \<      R  ,   \Diamond _{k}\partial _{p}\Phi _p  \> = \<     O(R^2),  \Diamond _{k}\Phi _p  \> ,\label{eq:mod eq}
  \end{align}
which need to be coupled with the following equation on $R$, obtained applying  $ P (p)$ to \eqref{eq:eq ansatz},
\begin{align}
  &    \im   \dot R  - {\mathcal H}_p  R     =   \( \dot \tau   -\lambda   (p) \) \cdot P (p) \sigma _3 \Diamond   R + \im \dot p\ P (p)\partial _p P (p)   R
 +    O(R^2) .\label{eq:eq for R}
  \end{align}
Equation \eqref{eq:eq for R} resembles a vectorial--like NLS.  Soffer and Weinstein in  \cite{SW1,SW2},  for a somewhat simpler system, had the idea to use
the dispersive properties of the linearized equation $ \im   \dot R  - {\mathcal H}_p  R     = 0$. Instrumental where advances in the dispersion  theory
of Schr\" odinger due to Journ\' e, Soffer and Sogge \cite{JSS}.
  Buslaev and Perelman in  \cite{BP1} for  dimension 1
  extended these results to the operator ${\mathcal H}_p $
  (the analysis in \cite{BP1} can be supplemented by material in \cite{KS}) introduced the idea of proving dispersion to 0 of $R$ by exploiting the dispersive properties
of  the group $e^{\im t {\mathcal H}_p}$.  Specifically, Buslaev and Perelman in  \cite{BP1}  prove   the following result.

\begin{theorem}[   Main Theorem in \cite{BP1} ]\label{thm:BP1}  In the $d=1$ dimension, suppose, in addition to the hypotheses in Theorem \ref{thm:class orb stab} that
     for  $\omega \in\mathcal{O}$ both the two  conditions listed below  are satisfied:
\begin{itemize} \item[(1)]   $0$ is the only eigenvalue of ${\mathcal H}_\omega$;
      \item[(H3)]    the points $\pm \omega$ are not resonances of  for ${\mathcal H}_\omega$.
\end{itemize}
Then for any
  $\omega _0\in \mathcal{O}$ there exists $\epsilon _0>0$  and $C_0>0$ s.t.
for
\begin{equation*}
      \|  \< x \> ^2 ( u _0 -
\phi _{\omega _0})   \| _{L  ^{2}} +  \|  \partial _x ( u _0 -
\phi _{\omega _0} )   \| _{L^2}<
\epsilon _0
\end{equation*}
we have
\begin{align*} &
     \| R(t) \| _{L^{\infty }} <C_0 (1+|t|)^{-1 /2}
\epsilon _0\ , \\& |  \dot \tau (t)  -\lambda   (p(t)) |+|\dot p(t)  |<C _0(1+|t|)^{-3 }\ \epsilon _0^2 \ .
\end{align*}
\end{theorem}
\noindent In \cite{Cu1} there is a version of the above result  for dimension $d\ge 3$, here quoted as is stated in  \cite{Cu1}.
\begin{theorem}\label{thm:Cu1}  For $d\ge 3$, assume the   hypotheses  of Theorem \ref{thm:BP1}. Then for any
  $\omega _0\in \mathcal{O}$ there exist  $\epsilon _0>0$  and $C_0>0$ s.t.,
for
\begin{align} &\label{eq:Cu1 in val1}
      \|  \< x \> ( u _0 -
\phi _{\omega _0})   \| _{H ^{2d+2[d/2]+2}}<
\epsilon _0  \\&    \|   u _0 -
\phi _{\omega _0}    \| _{H^{2d+2[d/2]+3} \cap W^{d+2[d/2]+2 ,1}}<
\epsilon _0 ,\label{eq:Cu1 in val2}
\end{align}
we have
\begin{align} &
     \| R(t) \| _{W^{d+[d/2]+1, \infty }} <C_0 (1+|t|)^{-d /2} \label{eq:Cu1 cont1}
\epsilon _0\ , \\& |  \dot \tau (t)  -\lambda   (p(t)) |+|\dot p(t)  |<C _0(1+|t|)^{-d }\ \epsilon _0^2 \ . \label{eq:Cu1 cont2}
\end{align}

\end{theorem}

Like in \cite{SW1,SW2,BP2}, crucial to \cite{Cu1}   is information on the  dispersion of the associated linearized evolution $e^{\im  t\mathcal{H}_\omega }$.  In fact, \cite{Cu1}
 contains the following theorem for $d\ge 3$, based on  work  by Yajima   \cite{Y1,Y2,galtabiar}, see also Weder, \cite{weder1,weder2},  which was inspired by Journ\' e, Soffer and Sogge \cite{JSS}. The case $d=2$ is in \cite{cuccagnatarullidim2}.

\begin{theorem}\label{thm:wave op}  For $d\ge2$, for any $\omega \in \mathcal{O}$ if, under the hypotheses (1) and (2)  of Theorem \ref{thm:BP1},
we set $L^q(\omega )= L^q(\R ^d, \C ) \cap N^{\perp}_g (\mathcal{L}_\omega ^*)$, then strong limits
\begin{equation}\label{eq:wave op1}
       W(\omega ) := s-\lim _{t\to +\infty }   e ^{  \im t \mathcal{H}_\omega}    e ^{ \im t\sigma _3 ( \Delta -\omega )} \text{ and } Z(\omega ) := s-\lim _{t\to +\infty }      e ^{\im t\sigma _3 ( -\Delta +\omega ) }e ^{ -  \im t \mathcal{H}_\omega}
\end{equation}
define    isomorphisms $L^q(\R ^d, \C ) \xrightarrow{  W(\omega ) }  L^q(\omega ) $    and  $L^q(\omega ) \xrightarrow{  Z(\omega ) }   L^q(\R ^d, \C ) $    for any $q\in [1, \infty ]$
and   yield   isomorphisms also between to Sobolev spaces $W ^{s,q}(\R ^d, \C ) $  and $    W ^{s,q}(\R ^d, \C ) \cap N^{\perp}_g (\mathcal{H}_\omega ^*) $ for any $s\in \R$.  Furthermore,
the   norms  of the operators are upper semicontinuous in  $\omega$.
\end{theorem}

\begin{remark} \label{rem:mistake1} Unfortunately, in \cite{Cu1} the proof of the case $q=2$, specifically  \cite[Corollary 3.2]{Cu1}, is wrong.  However, the correct proof  is a rather direct consequence of classical arguments by Kato \cite{kato}, and   is in \cite  {CPV}.

\end{remark}

The proof of Theorem \ref{thm:Cu1} involves applying ${\mathcal H}_p^j$  to \eqref{eq:eq for R} and then  applying to it $\< \cdot , \sigma _1\sigma _3\mathcal{H} _p {\mathcal H}_p^jR\>$
to get
\begin{align*}
  &      2^{-1} \partial _t \sum _{j=0}^{m}  \< {\mathcal H}_p^jR , \sigma _1\sigma _3\mathcal{H} _p {\mathcal H}_p^jR\>  = \sum _{j=0}^{m} \< {\mathcal H}_p^j\text{ r.h.s. of \eqref{eq:eq for R}}  , \sigma _1\sigma _3\mathcal{H} _p {\mathcal H}_p^jR\> ,
  \end{align*}
for a sufficiently large $m$. The summation in the l.h.s. is equivalent to $\| R \| _{H ^{2m+1}}^2 $. This is proved by induction from
$ \<  R ,  \sigma _1\sigma _3   R\>  \sim \| R \| _{H ^{ 1}}^2$ for $R\in N^{\perp}_g (\mathcal{L}_p ^*) $, which is true under the hypotheses of Theorem \ref{thm:Cu1}, see \cite{W2,CPV}.
By standard argument
\begin{equation}
      \| R (t) \| _{H ^{2m+1}}^2 \le C \| R (0) \| _{H ^{2m+1}}^2   e^{ \int _{0} ^t   \(    |  \dot \tau (t')  -\lambda   (p(t')) |+|\dot p(t')  |   +  \| R (t') \| _{W ^{m,\infty }}  \)  dt'}   . \label{eq:Cu1 gronw}
\end{equation}
This needs to be used in conjunction with estimates of the terms in the exponential. To this effect, we need to use {Theorem} \ref{thm:wave op}. Quite problematic  is   the   term  $ \( \dot \tau   -\lambda   (p) \) \cdot P (p) \sigma _3 \Diamond   R$   in \eqref{eq:eq for R}, as we will see also in Sect. \ref{subsec:dispersion f}.  Here we sketch the discussion in \cite{Cu1}, which comes from
  \cite{BP1}.  In an interval $[0,T]$,  we consider  $\dot \tau _1 = \lambda   (p(T)) $  with $\tau _1(0)=\tau (0)$ and define $R_1$ by $R_1 = e^{\im \sigma _3 (\tau -\tau _1)\cdot \Diamond} R$. Then, elementary computations
  yield
\begin{align}&
             P(p(T))  R  _1(t)  = e^ {\im t\mathcal H _{p(T)}}  P(p(T)) R _1 (0) \label{eq:Cu1 PpT}\\&  -  \im \int _{0} ^t  P(p(T)) e^ {\im (t-t')\mathcal H _{p(T)} } \( \dot \tau   -\lambda   (p) \)    \(   e^{\im \sigma _3 (\tau -\tau _1)\cdot \Diamond }  P(p) e^{ -\im \sigma _3(\tau -\tau _1)\cdot \Diamond } - 1    \)     \cdot   \sigma _3 \Diamond   R dt' \nonumber\\& -\im    \int _{0} ^t  P(p(T)) e^ {\im (t-t')\mathcal H _{p(T)} }\(   e^{\im \sigma _3 (\tau -\tau _1)\cdot \Diamond } V_p  e^{ -\im \sigma _3 (\tau -\tau _1)\cdot \Diamond } -V _{p(T)}  )       \) R_1dt'
 +  ...\nonumber
  \end{align}
Notice that, assuming \eqref{eq:Cu1 cont2},
\begin{align*}
              |\tau  -\tau _1|\le \int _0^t |\lambda (p) - \lambda(p(T))| dt' \lesssim \int _0^t  dt'  \int _{t'}^{T}|\dot p|  ds\lesssim   \epsilon _0^2 \int _0^t  dt'  \int _{t'}^{T}\< s \> ^{-2d}  ds \lesssim \epsilon _0^2.
  \end{align*}
This implies, by
\begin{align*}
             P(p(T))  R  _1(t) -R_1=  \left [  P(p(T)),   e^{ \im \sigma _3(\tau -\tau _1 )\cdot \Diamond}   \right ] R +  e^{\im \sigma _3 (\tau  -\tau  _1 ) \cdot \Diamond } \(  P(p(T))- P(p )     \) R,
  \end{align*}
that, assuming \eqref{eq:Cu1 cont1},  then the r.h.s. in the last equation is small, and so  $ P(p(T))  R  _1(t) \sim R  _1(t)$. In turn, by \eqref{thm:wave op}, which implies
\begin{equation*}
   \|  e^ {\im t \mathcal H _{p(T)}}  P(p(T)) R _1 (0) \| _{L^\infty} \lesssim \< t\> ^{-\frac{d}{2}}  \|    R _1 (0) \| _{L^1 \cap L^2} \text{ where }\< t\> : = \sqrt{1+|t|^2}
\end{equation*}
 we get, also from $\|    R _1 (0) \| _{L^1 \cap L^2} \lesssim \epsilon _0$ and from $ \< t\> ^{-\frac{d}{2}} * \< t\> ^{-\frac{d}{2}} \lesssim  \< t\> ^{-\frac{d}{2}} $,
\begin{align*}
          \|   P(p(T))  R  _1(t) \| _{L^\infty} \lesssim \< t\> ^{-\frac{d}{2}}  \epsilon _0+ \< t\> ^{-\frac{d}{2}}  \epsilon _0
  ^2.
  \end{align*}
This because the   terms in the last two lines of  \eqref{eq:Cu1 PpT}, being nonlinear, are smaller than the 1st term in the r.h.s. of   \eqref{eq:Cu1 PpT}.
Taking derivatives,  one gets  back \eqref{eq:Cu1 cont1}. Inserting this in the modulation equations \eqref{eq:mod eq}, one proves \eqref{eq:Cu1 cont2}.
This of course is just a caricature, but the rigorous argument is similar, assuming the \eqref{eq:Cu1 cont1}--\eqref{eq:Cu1 cont2} with some large  constant
$C_0$ and then proving, by taking $\epsilon _0$ sufficiently small, that the constant can be taken to be similar to
\begin{equation*}
    \sup _{\omega \in  \mathcal{O}  , \ t \ge 0}  \|  \< t \> ^{d/2}    e ^{  \im t \mathcal H_\omega} P(p(\omega ,0) ) \| _{L^1_x\cap L^2_x\to L^\infty _x} <\infty .
\end{equation*}
 Notice that \eqref{eq:Cu1 in val1}  is unnecessary (and is due to a non optimal choice in \cite{BP1,Cu1} of the modulation).

  \qed

Theorems   \ref{thm:BP1} and \ref{thm:Cu1}  are  based  in a significant way    on the fact that  $0$ is the only eigenvalue of ${\mathcal H}_\omega$.

 \begin{lemma} \label{thm:drop_(1)}   If we drop the hypothesis (1) in Theorems \ref{thm:BP1} and \ref{thm:Cu1}, then the conclusions of these theorems are false.

\end{lemma}
 \textit{Proof (sketch).}   Buslaev and Perelman showed, in \cite{BP2}, see also \cite{BS},
that the continuous spectrum component of $R(t)$ decays    slowly if ${\mathcal H}_\omega$ has exactly one eigenvalue $ \mathbf{e} (\omega ) \in (0,   \omega )$.
 We give a sketch
of this, assuming for simplicity that $u\in C^0(\R, H^1_{\rad}(\R ^d))$, thus excluding translations.
  Let us
suppose that $N  \lambda (\omega ) <\omega <(N+1)  \mathbf{e}  (\omega )$, for an $N\in \N$,  and let $\ker ({\mathcal H}_\omega -  \mathbf{e}  (\omega ))$ be generated by an appropriately normalized $\xi _\omega$.
Then, using the symmetry of $\sigma ({\mathcal H}_\omega)$ with respect to the coordinate  axes and the fact that  $\sigma ({\mathcal H}_\omega)=\sigma ({\mathcal H}_\omega ^*)$,
\begin{align}&
     r(t ) =   z  (t)  \xi    _{ \omega (t) }+  \overline{z}  (t) \sigma _1 { \xi }   _{ \omega (t) }         +f(t ) \text{  with } \label{eq:split_1} \\& f(t ) \in  \(   N_g ( {\mathcal H}_{ \omega (t) } )     \oplus  \ker ({\mathcal H}_{ \omega (t) } ^* -   \mathbf{e}  ({ \omega (t) })) \oplus \ker ({\mathcal H}_{ \omega (t) }^*+   \mathbf{e}  ({ \omega (t) } )) \) ^{\perp}.\nonumber
\end{align}
Then, in \cite{BP2} for the case $N=1$ and in \cite{zhousigal}  for generic $N$, it is shown that after a normal forms argument,
  for $P(|z|^2)$ real valued we have
\begin{equation} \label{eq:sistem_1}  \begin{aligned} &
{\rm i} \dot z -\mathbf{e}  z = P(|z|^2)z+  \overline{z}^N \langle f ,\sigma _1  G
( \omega )   \rangle _{L^2_x}+\cdots
\\& {\rm i}  \dot f -\mathcal{H}_\omega f = z^{N+1}
     M  ( \omega )+ \cdots
\end{aligned}
  \end{equation}
If,  in the equation for $z$, we substitute  $f$ with  $\displaystyle   -z^{N+1}\lim _{\varepsilon\to 0^+} R_{\mathcal{H}_\omega} \( (N +1)
 \mathbf{e}  (\omega ) +\im \varepsilon \) M  $, where the latter exists, see Proposition \ref{prop:resolvent} later,
 and use formula
 \begin{align}\label{eq:lim_abs_princ}
  \lim _{\varepsilon\to 0^+}  R_{\mathcal{H}_\omega}^{+}(\kappa +\im \varepsilon  ) = P.V. \frac{1}{\mathcal{H}_\omega -\kappa} - \im \pi \delta \( \mathcal{H}_\omega  - \kappa  \) \text{  for $\kappa \in \R $},
 \end{align}
which can be understood of the theory of distorted plane waves, but which we will not discuss in any detail,    then the equation of $z$
  becomes
\begin{equation}  \begin{aligned}      {\rm i}  \dot z-\mathbf{e}  (\omega ) z&= P(|z|^2)z- |z|^{2N }z\langle P.V. \frac{1}{\mathcal{H}_\omega -(N +1)
  \mathbf{e}  (\omega )} M , \sigma _1 G   \rangle _{L^2_x}  \\&{-{\rm i}|z|^{2N }z
\pi \langle    \delta ({\mathcal{H}_\omega}-(N +1)  \mathbf{e}  (\omega ) )M ,\sigma _1  G   \rangle} .\end{aligned}\nonumber
  \end{equation}
  Multiplying by  $\overline{z}$ and taking imaginary part, it can be shown that
  \begin{equation} \begin{aligned} \Rightarrow \frac{d}{dt}|z|^2&= - |z|^{2N+2 }  \Gamma (\omega) \text{   where }  \Gamma (\omega):= 2\pi   \langle    \delta ({\mathcal{H}_\omega}-(N +1)\, \mathbf{e}  (\omega ) )M ,\sigma _1  G   \rangle .  \end{aligned}
\nonumber
  \end{equation}
  Notice that, using an appropriate distorted Fourier transform associated to $\mathcal{H}_\omega$, we have
  \begin{equation}\label{eq:formula_gamma}
    \Gamma (\omega )\sim  \int _{ |\xi| = \sqrt{(N +1) \mathbf{e}  (\omega ) -\omega }}\< \widehat{M}(\xi
) , \sigma _1   {\widehat{G} (\xi )}\> _{\C ^2} dS.
  \end{equation}
  In the case $N=1$,  Buslaev and Perelman \cite{BP2} are able to show that the integral is   nonnegative.
  Zhou and Sigal  \cite{zhousigal} develop  rigorously the argument and assume  that $ \Gamma (\omega )>0$ to prove their own version of
  Theorem \ref{thm:main rough} for $N>1$ in the case of a single   $ \mathbf{e}  (\omega ) $. In   \cite{cuccagnamizumachi} it is shown that
  $  \Gamma (\omega )<0$  is incompatible with orbital stability. This means that, if there is a single $ \mathbf{e}   (\omega ) $,  if
  hypothesis (H1) and (H2) hold (they imply orbital stability), then, in the presumably generic case  $  \Gamma (\omega )\neq 0$, we need   to have
   $  \Gamma (\omega )> 0$. Notice that, if we assume $\Gamma (\omega) = \Gamma$  constant, then
   \begin{equation} \label{eq:formula z_t_1} \begin{aligned}  &|z(t)|  =  |z(0)| \( 1+ N\Gamma |z(0)| ^{2N} t\) ^ {-\frac{1}{2N}}. \end{aligned}
  \end{equation}
  The above discussion is purely heuristic, but indicative of the arguments in \cite{BP2,SW3,TY1,TY2,TY3,BS,Cu2,SW4,zhousigal,cuccagnamizumachi}.
 Notice that  by \eqref{eq:formula z_t_1}, eventually $|z(t)|\sim t  ^ {-\frac{1}{2N}} $ as $t\to +\infty$. In fact, since $|z(0)|$ is small,  $|z(t)|$ remains almost
   constant in the time interval $[0,|z(0)| ^{-2N} ]$. Because of the forcing term $z^{N+1}
     M $
 in \eqref{eq:sistem_1}, also $f$ cannot be counted to disperse  for a long time.  These arguments    show that the decay of $R(t)=P(p(t)r(t)$ in Theorem \ref{thm:Cu1}, in general cannot be
 expected to be true.

 \qed

The discussion in {Lemma} \ref{thm:drop_(1)} indicates the relevance of the eigenvalues of $\mathcal{H}_{\omega}$ in the analysis of the problem. In principle, eigenvalues of
$\mathcal{H}_{\omega}$ could lead to invariant tori  near the solitons, which would prevent the result in Theorem \ref{thm:main rough}. In fact we will discuss the fact that there
are no invariant tori, and this thanks to the a mechanism related to the fact that $ \Gamma (\omega )>0$  in \eqref{eq:formula_gamma}. The reader might wonder why we should have $ \Gamma (\omega )\ge 0$.
Heuristically this should be related to the fact that our NLS is Hamiltonian. If the coordinates $(z,f)$ in \eqref{eq:sistem_1}  were normal, we could expect   \eqref{eq:sistem_1} to be of the form
\begin{equation} \label{ham} \begin{aligned} &
{\rm i} \dot z = \partial   _{\overline{z}} E \, , \quad {\rm i}  \dot  f = \nabla
 _{\overline{f}} E.
\end{aligned}
  \end{equation}
  Then   by the Schwartz lemma, at  $z=0$ and $f=0$ we would get $(N+1)!M=\partial ^{N+1}_z\nabla _{\overline{f}}E
  =\overline{\partial ^{N }_{\overline{z}}\nabla _{f}\partial _{\overline{z}}E}=N! \sigma _1 \overline{{G}}$.  But   $M\sim \sigma _1 \overline{{G}}$ would imply $\widehat{M}\sim \sigma _1 \overline{\widehat{{G}}}$,
  yielding the $ \Gamma (\omega )\ge 0$.

\section{ The case when   $\mathcal{{H}}_{\omega }$   has positive eigenvalues  }
\label{sec:eigenvalues}
We will assume that $\Pi (\widetilde{u}_0) = \Pi (\widetilde{\phi} _{\omega ^1})=p^1$. This can be obtained using appropriate boosts.
We need some information on the spectrum of $\mathcal{H}_{\omega ^1}$. The following is elementary.

\begin{lemma}\label{lem:spectrum} For any $\omega \in \mathcal{O} $ following facts hold.
\begin{itemize}

\item[(1)] The spectrum $\sigma (\mathcal{{H}}_{\omega })$ is symmetric with respect to the coordinates axes. We have $\sigma (\mathcal{{H}}_{\omega })=\sigma (\mathcal{{H}}_{\omega } ^*)$.

\item[(2)] Since $\phi _\omega $ is a ground state,   all the eigenvalues  of $ \mathcal{{H}}_{\omega }$,
except   possibly for a pair $\pm   \im  \mathbf{e}  $ with  $ \mathbf{e}   >0$, are
in $   \mathbb{R}$.

\item[(3)] If $  \im  \mathbf{e}   \in \sigma _p (
 \mathcal{{H}}_{\omega })$ with $ \mathbf{e}   >0$ then
$N_g( \mathcal{{H}}_{\omega } - \im \mathbf{e}  ) $, the corresponding
generalized eigenspace, has
  dimension 1.

 \item[(4)]  If   $z\neq 0$ is an eigenvalue, then   we have $N_g(\mathcal{{H}}_{\omega }-z)=
 \ker (\mathcal{{H}}_{\omega }-z).$

\end{itemize}

\end{lemma}

We assume that $\mathcal{{H}}_{\omega ^1}$ has no embedded solitons inside the essential spectrum.

\begin{itemize}

 \item[(H4)]  There are no eigenvalues in $\mathcal{{H}}_{\omega ^1}$ in $ \R \backslash (-\omega ^1, \omega ^1 ) $.

\end{itemize}

\begin{remark} \label{rem:embedd1} It is expected, but unproved yet,  that, since $\phi _\omega$ is a ground state,  always there are no eigenvalues in  $\mathcal{{H}}_{\omega }$ in $\sigma _e(\mathcal{{H}}_{\omega }) = \R \backslash [-\omega, \omega ] $,
that is,  no eigenvalues embedded in the "interior" of  the continuous spectrum. Obviously, to our knowledge no  embedded eigenvalues  have been detected numerically in the case of ground states. For more general solitary waves which are not ground states,
we expect that  embedded eigenvalues could exist,  but that they cannot have positive \textit{Krein signature}.  The signature of the eigenvalues of  $\mathcal{H}_{\omega }$  in $ \R \backslash \{ 0  \} $ is
always positive, in the case of ground states.
\end{remark}

\begin{remark} \label{rem:embedd11}  Even in the case they exist, the embedded eigenvalues are unstable, in the sense that, perturbing the equation, the  $\mathcal{{H}}_{\omega }$ of the new equation
will in general not have these eigenvalues.  Results of this type go back to Grillakis \cite{G1}, are also in Tsai and Yau \cite{TY3} and, as explained  \cite{CPV}, can better be viewed
in the classical framework of Howland \cite{how1,how2}
\end{remark}

We allow  $\mathcal{{H}}_{\omega ^1}$ to have a certain number of eigenvalues in the gap $(-\omega _1, \omega _1)$.

\begin{itemize}

 \item[(H5)]  There is an $m$ s.t.\   $\mathcal{H}_{\omega ^1}$ has  $m$ positive eigenvalues  $\mathbf{e}
_1 \le \mathbf{e}   _2 \le ...\le \mathbf{e}   _m $, where we repeat an eigenvalue a number of times equal to
its multiplicity.
We assume there are fixed integers $m_0=0< m_1<...<m_{l_0}=m$ such that
$\mathbf{e}   _j = \mathbf{e}   _i $ exactly for $i$ and $j$
both in $(m_l, m_{l+1}]$ for some $l\le l_0$. In this case $\dim
\ker (\mathcal{H}_{\omega ^1 }-\mathbf{e}   _j  )=m_{l+1}-m_l$.
 We assume there exist $N_j\in \mathbb{N}$ such that
$0<N_j\mathbf{e}   _j < \omega _1< (N_j+1)\mathbf{e}   _j $ with
$N_j\ge 1$. We set $N=N_1$.

\end{itemize}

\begin{remark} \label{rem:mult_liter}   The literature considered for more than a decade only the case when $m=1$,   except
for \cite{T}, where however  $N_j=1$  for all $j$. These are very restrictive conditions.
   Only \cite{bambusicuccagna,Cuccagna11CMP} started to consider fairly general situations.
 \end{remark}

\begin{remark} \label{rem:mult}  Hence we allow the eigenvalues to have finite multiplicity. The number $(N_j+1)\in \N$ is the smallest such that the corresponding multiple of
of $\mathbf{e}   _j $ is in $\sigma _c (\mathcal{H}_{{\omega ^1}})$.
 \end{remark}

\begin{remark} \label{rem:numerical1}  We give here a partial list of papers which have explored   the spectrum of
operators such as $\mathcal{{H}}_{{\omega ^1}}$.  Chang,  Gustafson, Nakanish and   Tsai  \cite{Chang} explore in great detail and mostly numerically   the spectrum of $\mathcal{{H}}_{{\omega ^1}}$ in the case
of $\beta (|u|^2)=-|u|^{p-1} $.  Their computations in dimensions $d=1,2,3$ show the presence of many eigenvalues for $p\to 1 ^{+}$  and of just two real nonzero  eigenvalues for $p\to \( 1+4/d\) ^{-}$
which reach 0   at $p= 1+4/d$ and bifurcate into    two imaginary      eigenvalues    for $p> 1+4/d $.
 \end{remark}

 \begin{remark} \label{rem:numerical2}    Buslaev and Grikurov
\cite {BG01MCS}  and    Marzuola,   Raynor, and  Simpson   \cite{MRS10JNS} study numerically situations when $q(\omega )$,  the function in (H2) Theorem \ref{thm:class orb stab}, has a minimum $\omega _*$.
Then $\mathcal{H}_{\omega }$ has two imaginary eigenvalues for $\omega <\omega _*$ which converge to 0 as $\omega \to \omega _* ^-$ and bifurcate into two positive eigenvalues for  $\omega >\omega _*$.
This is explained analytically in  Comech  and Pelinovsky  \cite{CP03CPAM}.  Interesting oscillating patterns are described numerically in   \cite {BG01MCS,MRS10JNS},
with interesting conjectures, which are discussed analytically, but inconclusively, in  \cite{CM19JDE}, where the problem is shown to be similar to   that of a soliton constrained in a potential.
 \end{remark}

\begin{remark} \label{rem:numerical3}  The spectrum of $\mathcal{H}_{{\omega ^1}}$ for
the equations with $\beta (|u|^2)=-|u|^{p-1} $ and $p>1+4/d$, in particular the case $d=3$ and $p=3$ have been studied in considerable detail.
The case $\beta (|u|^2)=-|u|^{2} $, $d=3$, is considered in Schlag \cite{Schlag}, where it is shown that  $\mathcal{H}_{{\omega ^1}}$ has no eigenvalues other than
0 in $[-{\omega ^1}, {\omega ^1}]$ if the operator $L_{+{\omega ^1}}$ in Theorem \ref{thm:class orb stab} and the operator $L_{-\omega _1}:= -\Delta -{\omega ^1}+ \beta  (\phi _{\omega ^1}^2)  $
don't have eigenvalues in  $(0, {\omega ^1}]$. This information on $ L_{\pm {\omega ^1}}$
is verified numerically for the cubic NLS  with $d=3$    in  Demanet and Schlag \cite{Demanet} and proved rigorously  in Costin, Huang and Schlag \cite{Costin}. In Marzuola and Simpson  \cite{Marzuola}, for the cubic NLS  with $d=3$ it is proved numerically absence of nonzero real eigenvalues.
 Further cases  of  computer assisted proofs of absence of   nonzero real eigenvalues for mass supercritical   NLS   with   $\beta (|u|^2)=-|u|^{p-1} $   are considered in  Asad and Simpson  \cite{asad}.

 \end{remark}

We assume that the eigenvalues in (H5) satisfy the following non resonance condition.
\begin{itemize}

\item [(H6)]  If $\mathbf{e}   _{j_1}<...<\mathbf{e}   _{j_k}$ are $k$ distinct
  $\mathbf{e}  $'s, and $\mu\in \Z^k$ satisfies
  $|\mu|:=|\mu_1|+...+|\mu_m| \leq 2N_1+3$, then we have
$$
\mu _1\mathbf{e}   _{j_1}+\dots +\mu _k\mathbf{e}   _{j_k}=0 \iff \mu=0\ .
$$

\end{itemize}

\begin{remark} \label{rem:lin_indip}   A more restrictive   formulation  would be to say that the eigenvalues are linearly independent in $\Z$. That would be a more stringent condition
than necessary.
 \end{remark}
Another hypothesis it the following.

\begin{itemize}

\item [(H7)]  There is no multi index $\mu \in \mathbb{Z}^{m}$
with $|\mu| \leq 2N_1+3$ such that $\mu \cdot
\overrightarrow{\mathbf{e}} =\omega ^1$   (where $\overrightarrow{{\mathbf{e}}}:=(\mathbf{e}_1,..., \mathbf{e}_m))$.

\end{itemize}

\begin{remark} \label{rem:hyp_omega}   Notice that in \cite{Cuccagna11CMP} and in some of the subsequent papers, the hypotheses are more restrictive, because
it is assumed that the multiplicities of the eigenvalues are constant, and the hypotheses (H4)--(H7) are assumed for all $\omega$.  The hypotheses stated here,
which require  (H4)--(H7)   just for $\omega ^1$,
come from \cite{Comech}.
 \end{remark}

We need to record the following version of Theorem \ref{thm:wave op}, which has essentially the same proof, see \cite{CPV} on how to deal with the eigenvalues.
\begin{theorem}\label{thm:wave op1}  Let
\begin{align}  \label{eq:spectraldecomp} & X_c (\omega ^1 ):=
\left\{N_g(\mathcal{H}_{\omega ^1 } ^\ast)\oplus \big (\oplus _{ \mathbf{e}   \in \sigma
_p\backslash \{ 0\}}   \ker (\mathcal{H}_{\omega ^1 } ^*- \mathbf{e}
 ) \big)\right\} ^\perp ,
\end{align}
where we can take $X_c  (\omega ^1 )\subset \mathcal{S}^{\prime}(\R ^d, \C^ 2)$, in the space of tempered distributions. Then the statement of Theorem \ref{thm:wave op} continues to be true
 for  $L^q(\omega ):= L^q(\R ^d, \C ^2) \cap X_c  (\omega ^1 )$.  In particular the wave operators    yield   isomorphisms also between to Sobolev spaces $W ^{s,q}(\R ^d, \C ^2 ) $  and $    W ^{s,q}(\R ^d, \C ^2 ) \cap X_c  (\omega ^1 )$ for any $s\in \R$.
\end{theorem}
Notice that
\begin{align} 	\label{eq:spectdec-1} & N_g^\perp
(\mathcal{H}_{\omega ^1}^{\ast})   =  \big (\oplus _{\mathbf{e}  \in \sigma _p\backslash \{ 0\}}
\ker (\mathcal{H}_{\omega ^1}  -\mathbf{e}
 ) \big) \oplus (L^2(\R ^d, \C ^2) \cap X_c  (\omega ^1 ))
\end{align}
and so  the spectral decomposition
\begin{align} 	
\label{eq:spectdec0}& L^2(\R ^d, \C ^2) = N_g(\mathcal{H}_{\omega ^1})\oplus  \big (\oplus _{\mathbf{e}  \in \sigma _p\backslash \{ 0\}}
\ker (\mathcal{H}_{\omega ^1}  -\mathbf{e}
 ) \big) \oplus (L^2(\R ^d, \C ^2) \cap X_c  (\omega ^1 )).
\end{align}
 We have the following useful result which, among other things, insures that $\mathcal{H}_{\omega ^1}$  satisfies the limiting absorption principle.
\begin{proposition}
  \label{prop:resolvent}  There exists $\tau _d>0$ s.t.\ for $\tau \ge \tau _d$ the following hold.

{\item {(1)}} There exists $C=C(\tau ,\omega )$, upper
semicontinuous in $\omega $ such that  for any $\varepsilon \neq 0$,
$$\| R_{\mathcal{H}_{\omega  } }(\lambda +\im \varepsilon )P_c(\mathcal{H}_\omega )
u\| _{L^2_\lambda L^{2,-\tau }_x}\le  C \| u\| _{L^2}.  $$
where
\begin{equation*}
     L^{2, \tau }_x = L^{2,-\tau }(\R ^d, \C^2) =\{ f \in \mathcal{S}' (\R ^d, \C ^2) : \|  (1+|x|^2)^{-\frac{\tau }{2}}f\| _{L^2}<\infty \} .
\end{equation*}

{\item{(2)}} For any $u\in L^{2, \tau }_x $ the following limits exist:
$$ \lim _{\epsilon \searrow 0}R_{\mathcal{H}_\omega }(\lambda \pm \im \varepsilon )
u= R_{\mathcal{H}_\omega }^\pm  (\lambda ) u  \text{ in $C^0(\sigma
_e(\mathcal{H}_\omega ),L^{2, -\tau }_x)$}.$$

{\item {(3)}} There exists $C=C(\tau ,\omega )$, upper
semicontinuous in $\omega $ such that
$$
  \|   R_{\mathcal{H}_\omega }^\pm  (\lambda
)P_c(\mathcal{H}_\omega )   \| _{B( L^{2,\tau }_x, L^{2,-\tau }_x)} < C
\langle \lambda \rangle ^{-\frac{1}{2}} .$$

{\item {(4)}} Given any
$u\in L^{2, \tau }_x $ for the projection on the $L^2(\R ^d, \C ^2) \cap X_c  (\omega ^1 )$   term in \eqref{eq:spectdec0}    we have
$$P_c(\mathcal{H}_\omega )u=\frac{1}{2\pi \im }\int _{\sigma _e(\mathcal{H}_\omega )}
(R_{\mathcal{H}_\omega }^{+}(\lambda  )-R_{\mathcal{H}_\omega }^{-}(\lambda  ))  u\,
d\lambda .$$
\end{proposition}

\section{ Idea of the  proof of   {Theorem} \ref{thm:main rough}}
\label{sec:eigenvalues}

In   this section, we will proceed to show heuristically how  to prove {Theorem} \ref{thm:main rough}.

First of all, by some small boosts  we reduce to the case when $\Pi _a(u_0)=0$ for all $a=1,...,d$. We notice that the we can write $R=P(p)r$ with $r\in N^{\perp}_g (\mathcal{H} _{\omega_{1}}  ^*)$.
Now we have
\begin{equation} \label{eq:SystE1} \begin{aligned} &
\dot p = \{ p, E \} \,, \quad \dot \tau = \{ \tau, E \} \,,  \quad \dot r= \{ r, E
\}. \end{aligned}
\end{equation}
We can substitute the coordinates $p$ with the coordinates $\Pi$. In the new coordinates, the system becomes
\begin{align} &
\dot \Pi_j = 0 \ , \quad \dot \tau = \{ \tau, E \} \,, \nonumber \\&   \quad \dot r= \{ r, E
\}.   \label{eq:SystE2}\end{align}
Notice that in the   coordinates  $(\Pi , \tau , r)$, as well as in the system of coordinates  $(p , \tau , r)$,
we have $\partial _\tau E=0$. Then we have a reduction of the system to $\dot r= \{ r, E
\}.$

We also choose $p^0\in \mathcal{P}$ so that
\begin{equation}\label{eq:coordinate-3} \text{$ \Pi ( u_0)= p^0 $  }.
\end{equation}
Notice that if we consider the equations $\Pi  =   p^{0} $, they define a submanifold in $H^1(\R ^d, \widetilde{\C} )$
in a neighborhood of $\{   e^{\im \tau \cdot \Diamond } \Phi _{p^1} : \tau  \in \R ^{d+1} \}$
This set is parametrized by $(\tau , r)$. Taking the quotient by the group $e^{\im \tau \cdot \Diamond }$,
we obtain a manifold, parametrized by $r$. This manifold inherits in a natural way
a symplectic structure, which is inherited from the $\Omega$ defined in \eqref{eq:SymplecticForm}.

Now we split according to \eqref {eq:spectdec-1}
 \begin{equation}
 \label{eq:decomp2}
 r (x) =\sum_{l=1,..., {{n}}}z_l \xi_l (x) +
\sum_{l=1,..., {{n}}}\overline{z}_l  \sigma _1{{\xi }}_l (x)
+ f (x), \quad f \in X_c \text{ with $f= \overline{{f}} $ },
\end{equation}
We will assume that it is possible to change the $(z,f)$ coordinate so that the symplectic form is
given by
\begin{equation}\label{eq:modomega0}
 \Omega= \im \sum_{l=1,...,n} \,dz_l\wedge d \overline{z}_l
+\im  \langle   \sigma_3 df, \sigma _1df\rangle .
\end{equation}
The correct version is just slightly more complex, see \cite[formula (7.11)]{Comech}, and contains some additional higher order terms that we skip in the following  heuristic  discussion.

\noindent It is critical, to consider an expansion of the energy in terms of these coordinates. It is important to expand
\begin{align*}  E &=  E(\Phi_{p}) + \< \nabla E(\Phi_{p}), \sigma _1 P_{p}r\> + 2^{-1}\< \nabla  ^2E(\Phi_{p}) P_{p}r, P_{p}r \>+...\\& = E(\Phi_{p})   + 2^{-1}\< \nabla  ^2E(\Phi_{p}) P_{p}r, P_{p}r \> +... \\&
\end{align*}
where we  exploit  $ \< \nabla E(\Phi_{p}), \sigma _1 P_{p}r\>  = \lambda  \cdot  \<  \Diamond \Phi_{p} , \sigma _1 P_{p}r\>  =0.$  Adding and subtracting $\lambda (p) \cdot  \Pi  (P_{p}r)$ in the r.h.s.,
we obtain
\begin{align*}  E  = E(\Phi_{p}) +\lambda (p) \cdot  \Pi  (P_{p}r)  + 2^{-1}\< \(\nabla  ^2E(\Phi_{p})-  \lambda  \cdot \Diamond  \) P_{p}r, P_{p}r \> +...
\end{align*}
Substituting $\Pi = p+  \Pi  (P_{p}r)$,   subtracting on both sides $E(\phi_{\omega ^0})$ ,   we get
\begin{align*}  E  -E(\phi_{\omega ^0}) & = E(\Phi_{p}) -\lambda   \cdot  p -\( E(\phi_{\omega ^0}) -\lambda _0  \cdot  p ^0\)   +(\lambda -\lambda _0)  \cdot  p^0     \\&+ 2^{-1}\< \(\nabla  ^2E(\Phi_{p})-  \lambda  \cdot \Diamond  \) P_{p}r, P_{p}r \> +... \\&= d(\omega )- d(\omega ^0)+(\omega -\omega ^0)q(\omega ^0) + 2^{-1} v^2 q(\omega ^0) + 2^{-1}\< \sigma _3 \mathcal{H}_{\omega ^1}  r,  r \>+ ...
\end{align*}
 Let us now substitute $r$ with the expansion  in \eqref{eq:decomp2}. Then, for $d(\omega ):= E(\phi_{\omega  })-\omega q(\omega )$,
\begin{align*}  E  -E(\phi_{\omega ^0}) & = E(\Phi_{p}) -\lambda   \cdot  p -\( E(\phi_{\omega ^0}) -\lambda _0  \cdot  p ^0\)   +(\lambda -\lambda _0)  \cdot  p^0     \\&+ 2^{-1}\< \(\nabla  ^2E(\Phi_{p})-  \lambda  \cdot \Diamond  \) P_{p}r, P_{p}r \> +... \\&= d(\omega )- d(\omega ^0)+(\omega -\omega ^0)q(\omega ^0) + 2^{-1} v^2 q(\omega ^0) + 2^{-1}\< \sigma _3 \mathcal{H}_{\omega ^1}  r,  r \>+ ...
\end{align*}
   Then we obtain an expression of the form
\begin{align*}   E  &= \psi (\Pi (f))  +E _{\text{discr}} + 2^{-1}\< \sigma _3 \mathcal{H}_{\omega ^1}  f, \sigma _1 f \> \\&  +  \sum _{ \substack{
 }} z^\mu \overline{z}^\nu a_{\mu\nu } (\Pi (f)) +   \sum _{ \substack{
 }} z^\mu \overline{z}^\nu  \<  \sigma _3 A_{\mu\nu } (\Pi (f)) ,\sigma _1f \>   + ...
\end{align*}
 where $E _{\text{discr}} :=\sum _j \mathbf{e} _j |z_j|^2$ and
where we sum over finitely many multi--indexes.  We   remark that $\overline{a}_{\mu\nu }=a_{\mu\nu }$  and $ \overline{A}_{\mu   \nu} =-\sigma _1{A}_{\nu   \mu}$, by the fact that $E$ is real valued.

\noindent Non resonant terms  of the form $z^\mu \overline{z}^\nu a_{\mu\nu }  $  for $(\mu -\nu )\cdot \overrightarrow{\mathbf{e}} \neq 0$can be eliminated by   considering
 appropriate canonical transformations given by $\phi ^{t} |_{t=1}$, using the flow of the Hamiltonian vector--field associated
to functions of the form $\chi =z^\mu \overline{z}^\nu b_{\mu\nu } $, with the coefficient unknown. Indeed, concisely,
\begin{align*} & E \circ  \phi ^{t} |_{t=1}=  E +\{ E _{\text{discr}} + 2^{-1}\< \sigma _3 \mathcal{H}_{\omega ^1}  f, \sigma _1 f \>,   z^\mu \overline{z}^\nu  \} b_{\mu\nu }+...\\& =  E +\{ E _{\text{discr}}   \} b_{\mu\nu }+.. =E +  \overrightarrow{\mathbf{e}}\cdot (\mu -\nu )z^\mu \overline{z}^\nu b_{\mu\nu }+...
\end{align*}
can be used to eliminate the non resonant  $z^\mu \overline{z}^\nu a_{\mu\nu }  $ term, just by solving $a_{\mu\nu } +\overrightarrow{\mathbf{e}}\cdot (\mu -\nu )  b_{\mu\nu } =0$.

\noindent Similarly, terms of the form $z^\mu \overline{z}^\nu  \<  \sigma _3 A_{\mu\nu }   ,\sigma _1f \>$ with $ |\overrightarrow{\mathbf{e}}\cdot (\mu -\nu )|<\omega ^{1}$  are \textit{non resonant}, and can be eliminated
similarly using  $\chi =z^\mu \overline{z}^\nu  \<  \sigma _3 B_{\mu\nu }   ,\sigma _1f \> $. Indeed, concisely,
\begin{align*} & E \circ  \phi ^{t} |_{t=1}=  E +\{ E _{\text{discr}} ,  z^\mu \overline{z}^\nu   \}  \<  \sigma _3 B_{\mu\nu }   ,\sigma _1f \> +  z^\mu \overline{z}^\nu    \{  2^{-1}\< \sigma _3 \mathcal{H}_{\omega ^1}  f, \sigma _1 f \> ,  \<  \sigma _3 B_{\mu\nu }   ,\sigma _1f \> \}  +...\\& =E +  z^\mu \overline{z}^\nu  \<  \sigma _3 \( \overrightarrow{\mathbf{e}}\cdot (\mu -\nu )  + \mathcal{H}_{\omega ^1}\)   B_{\mu\nu }  ,\sigma _1f \>+...
\end{align*}
so that the non resonant term can be canceled solving $\( \overrightarrow{\mathbf{e}}\cdot (\mu -\nu )  + \mathcal{H}_{\omega ^1}\)   B_{\mu\nu }  =A_{\mu\nu }  $.  Here, the fact that the coefficients
$a_{\mu\nu }$  and  $A_{\mu\nu }$, and so also $b_{\mu\nu }$  and  $B_{\mu\nu }$, depend on $\Pi (f)$ and are not constant, is not an obstacle for a rigorous implementation of the above ideas, because
$\Pi (f)$ remain constant, up to an error which is higher order and does not affect the computations.

 \noindent Eventually we find a a system of coordinates, where the significant terms are
\begin{align*}E =&\psi (\Pi (f)) +E _{\text{discr}} + 2^{-1}\< \sigma _3 \mathcal{H}_{\omega ^1}  f, \sigma _1 f \> + Z_0   + Z_1+... \text{ , where }\\& Z_0=  \sum _{ \substack{
 (\mu -\nu )\cdot \overrightarrow{\mathbf{e}}=0}} z^\mu \overline{z}^\nu a_{\mu\nu } (\Pi (f)) \\& Z_1:=\sum _{\overrightarrow{\mathbf{e}}
 \cdot \mu >\omega ^1}z^\mu   \< \sigma _3A_{\mu   0},\sigma _1f  \> + \sum _{ \overrightarrow{\mathbf{e}}
\cdot  \nu >\omega ^1}  \overline{z}^\nu \< \sigma _3A_{0   \nu},\sigma _1f  \>  .
\end{align*}
Here $E$ real valued, we have $\overline{f}= \sigma _1f$, and as a consequence
\begin{align*}  \overline{A}_{\mu   \nu} =-\sigma _1{A}_{\nu   \mu}.
\end{align*}
The system reads
\begin{align} &
   {\rm i}  \dot f  =  \mathcal{H}_{\omega  ^1}   f  + \sigma _3  \nabla _{\Pi (f)}E \cdot \Diamond f+ \sum  _{
  \mathbf{e}  \cdot  \alpha    > \omega ^1
  }   z ^\alpha
             A _{\alpha 0}    + \sum  _{
  \mathbf{e} \cdot  \beta    > \omega ^1
  }
 \overline{ {z }}^ { \beta }           A _{0 \beta} +...  \label{eq:eq for f}\\&    {\rm i}    \dot z _j
= \partial _{\overline{z}_j}  E=    \mathbf{e} _j z_j   +\partial _{\overline{z}_j} Z_0+    \sum  _{
  \overrightarrow{\mathbf{e}}  \cdot \nu   > \omega ^1
  }  \nu _j\frac{
 \overline{ {z }}^ { {\nu} } }{\overline{z}_j} \langle  A_{0 \nu } ,\sigma _3\sigma _1f \rangle  +...\label{eq:eq for f}
\end{align}
The crux of the proof consists in proving the following.

 \begin{proposition}\label{proposition:mainbounds} There is a fixed
$C _0>0$ such that for $\varepsilon _0>0$ sufficiently small, for
$\epsilon \in (0, \varepsilon _0)$ and for $|z(0)|+\| f (0)\| _{H^1}<\epsilon$,
then the following inequalities, for some $T>0$
\begin{align}
&   \|  f \| _{L^r_t( [0,T]  ,W^{ 1 ,p}_x)}\le
  2C _0\epsilon \text{ for all admissible pairs $(r,p)$}
  \label{Strichartzradiation}
\\& \| z ^\mu \| _{L^2_t( 0,T )}\le
  2C _0\epsilon \text{ for all multi indexes $\mu$
  with  $ \overrightarrow{\mathbf{e}} \cdot \mu >\omega _0 $} \label{L^2discrete}\\& \| z _j  \|
  _{W ^{1,\infty} _t  ( 0,T   )}\le
  2C_0 \epsilon \text{ for all   $j\in \{ 1, \dots , m\}$ }
  \label{L^inftydiscrete}
\end{align}
imply improved  inequalities obtained replacing $2C_0$  with $C_0$.

\end{proposition}

\subsection{Analysis of the equation of $f$}
\label{subsec:dispersion f}
If we had $ \sigma _3  \nabla _{\Pi (f)}E \cdot \Diamond f$=0, then we would have
\begin{align} &
 \|  f \| _{L^r_t( [0,T]  ,W^{ 1 ,p}_x)} \le c_0 \| f (0)\| _{H^1}  +   \sum  _{
  \overrightarrow{\mathbf{e}}  \cdot \mu   > \omega ^1
  }   \| z ^\mu \| _{L^2_t( 0,T )} + C(C_0) \epsilon ^2 ,   \label{eq:disp f1}
\end{align}
 implying that the key estimates are those on   $\| z ^\mu \| _{L^2_t( 0,T )}$ for $\overrightarrow{\mathbf{e}}\cdot \mu   > \omega ^1$.

 \noindent  This in fact is true, but nonetheless
\begin{align} &
  \sigma _3  \nabla _{\Pi (f)}E \cdot \Diamond f = \varpi _0\sigma _3f+ \im \sigma _3\overrightarrow{\varpi} \cdot \nabla f  \label{eq:int_fac_1}
\end{align}
is   nonzero. It helps that the hypotheses of Proposition  \ref{proposition:mainbounds} imply  $\| ( \varpi _0,\overrightarrow{\varpi} )\| _{L^\infty (0,T)} \le C(C_0)  \epsilon $.  However
the terms in \ref{eq:int_fac_1} are  not in  $  L^2([0,T], W^{1,\frac{2d}{d+2}}) +  L^1([0,T], H^{1} ) $, so that they cannot be incorporated with the terms on their   right in \eqref{eq:eq for f} .  Nor  they can be
eliminated easily by some integrating factor, this because    $\sigma _3$   and $\im \sigma _3 \partial _a$ for $a=1,...,d$ do not commute with
$ \mathcal{H}_{\omega  ^1}$.

Nonetheless, a form of integrating factor has been proved by Beceanu \cite{beceanu}, but only in dimensions $d\ge 3$. An alternative argument, attributed to Perelman,
is presented in Bambusi \cite[Appendix B]{Bambusi13CMPas}, but that too, based on Proposition 1.1 \cite{perelman3}, depends on dimensions $d\ge 3$.
A different argument, due to
Buslaev and Perelman \cite{BP2} is known in dimensions $d=1,2$, but only when $\overrightarrow{\varpi}=0$  in \eqref{eq:int_fac_1}, that is, when there is no translation in the problem.
This accounts for the fact that  {Theorem} \ref{thm:main rough} has not been proved in dimensions 1 and 2, as we already mentioned, except under hypotheses (like extra symmetries, or in the presence of a
potential) that break the translation invariance.   Notice that the gauge change argument sketched near   \eqref{eq:Cu1 PpT} depends   on the hypothesis  (1)   stated in Theorem \ref{thm:BP1} (absence of nonzero eigenvalues).

\begin{remark} \label{rem:new_dim1}  What is crucial, for the integrating factor argument, is that
\begin{equation*}
   \| \< x \> ^{-M(d)} e^{\im t \Delta } \< x \> ^{-M(d)} \| _{L^2\to L^2} \le C_d \< t \> ^{-\frac{d}{2}}
\end{equation*}
	 for a sufficiently large $M(d)$, with  $\frac{d}{2}>1$ for $d\ge 3$. In the cases $d=1,2$ the lack of integrability is an obstruction for the argument. The case $d=1$ in particular, can be phrased by stating that
$-\Delta$ is in $\R$ a non--generic Schr\"odinger operator, because of the fact that the point 0 is a resonance.

It is probably not coincidental that the proof of asymptotic stability of   kinks for the $\phi ^4$ model by  Kowalczyk,   Martel and  Mu\~{n}oz \cite{KMM2} is valid only in the case of odd solutions, that is, by imposing
a symmetry which allows to exclude  translation, and that the main difficulty at removing the  symmetry  is the   resonance at the threshold of the continuous spectrum of the linearization, see \cite[Remark 1.2]{KMM2}.
\end{remark}

 \subsection{Analysis of the equation of $z$}
\label{subsec:FGR z}

Recall that we defined   $E _{\text{discr}} :=\sum _j \mathbf{e} _j |z_j|^2$.  The idea of the proof consists, schematically,   in showing that
\begin{align} & \dot E _{\text{discr}} =\left \{ E _{\text{discr}}, E \right \}
\sim \left \{ E _{\text{discr}}, Z_1 \right \}  \lesssim - \sum _{\overrightarrow{ \mathbf{e}}
 \cdot \mu >\omega ^1}|z^\mu |^2    \label{eq:fermi1}
\end{align}
This will imply
\begin{align*} &   \sum _j \mathbf{e} _j |z_j(t)|^2+  \sum _{ \mathbf{e}
 \cdot \mu >\omega ^1} \int _0^t|z^\mu |^2  \lesssim \sum _j \mathbf{e} _j |z_j(0)|^2,
\end{align*}
yielding the crucial  bound that, in turn and thanks to \eqref{eq:disp f1}, yields Proposition  \ref{proposition:mainbounds}. We sketch a heuristic argument
for \eqref{eq:fermi1} The starting point consists in considering
 \begin{equation}
\begin{aligned} &
  f =- \sum   z ^\alpha
   R ^{+}_{ \mathcal{H}_{\omega ^1}  }(\overrightarrow{\mathbf{e}}  \cdot  \alpha    ) A _{\alpha 0} - \sum
 \overline{ {z }}^ { \beta  } R ^{+}_{ \mathcal{H}_{\omega ^1}  }(    -\overrightarrow{\mathbf{e}}  \cdot \beta ) A _{0 \beta} +g. \label{eq:fermi2}
\end{aligned}
\end{equation}
The effect of this change of variables is to show that $g$ satisfies an equation  where the terms $z ^\alpha  A _{\alpha 0}$  and $\overline{z} ^\alpha  A _{0\alpha  }$ have been canceled out.
While $g(t)\not \in H^1$, nonetheless, using \eqref{Strichartzradiation}--\eqref{L^inftydiscrete} it is possible to prove
\begin{align*}
&   \|  \< x \> ^{-N(d) } g \| _{L^2_t( [0,T]  ,H^1_x)}\lesssim \| f(0)\|  _{H^1}+|z(0)| +\epsilon ^2 \le C_0  \epsilon
\end{align*}
for an appropriate $N(d)>0$.  Here $c_0\sim 1$ is $c_0\ll C_0$ and consequently in the sequel we ignore
the term $g$.

\noindent Substituting  \eqref{eq:fermi1}  in \eqref{eq:eq for f}, ignoring the terms in $g$, which are smaller,  by elementary arguments we get to

\begin{equation*}
\begin{aligned} &
   {\rm i}   \dot z _j
= \mathbf{e} _j z_j -    \sum  _{\substack{
 \overrightarrow{\mathbf{e}}  \cdot  \nu  > \omega ^1
  \\  \overrightarrow{\mathbf{e}} \cdot  \alpha     > \omega ^1}}  \nu _j\frac{z ^\alpha
 \overline{ {z }}^ { {\nu} }    }{\overline{z}_j} \langle  A_{0 \nu } ,\sigma _3\sigma _1R ^{+}_{ \mathcal{H}_{\omega ^1}  }(\overrightarrow{\mathbf{e}}  \cdot  \alpha ) A _{\alpha 0 }   \rangle +  ...
\end{aligned}
\end{equation*}
Generically, when all  the eigenvalues  of $\mathcal{{H}}_{\omega _1}$ in  $(0, \omega _1)$  have multiplicity 1, and recalling   $A_{0 \alpha }=-\sigma _1 \overline{A}_{  \alpha 0}  $, this simplifies further
\begin{equation}
\begin{aligned} &
   {\rm i}    \dot z _j
=  \mathbf{e} _j z_j -   \sum  _{\substack{  \overrightarrow{\mathbf{e}} \cdot  \alpha     > \omega ^1}}  \alpha _j\frac{|z ^\alpha | ^2
     }{\overline{z}_j} \langle \overline{A} _{\alpha 0} ,\sigma _3 R ^{+}_{ \mathcal{H}_{\omega ^1}  }( \overrightarrow{\mathbf{e}}\cdot  \alpha ) A _{\alpha 0}   \rangle  +... .\label{eq:last eq z}
\end{aligned}
\end{equation}
Recalling now  {Theorem} \ref{thm:wave op},  for $A _{\alpha 0} =   W(\omega ^1)A _{\alpha 0}^{(1)}$, we have the following steps, already discussed
in the old survey \cite{cuccagna2004},
\begin{align*}
&     \langle \overline{A} _{\alpha 0} ,\sigma _3 R ^{+}_{ \mathcal{H}_{\omega ^1}  }( \overrightarrow{\mathbf{e}}\cdot  \alpha ) A _{\alpha 0}   \rangle  = \lim _{\varepsilon \to 0^+}\langle \overline{A} _{\alpha 0} ,\sigma _3 R _{ \mathcal{H}_{\omega ^1}  }( \overrightarrow{\mathbf{e}}\cdot  \alpha +\im \varepsilon) A _{\alpha 0}   \rangle \\& = \lim _{\varepsilon \to 0^+}\langle  \overline{W(\omega ^1)A _{\alpha 0}^{(1)}} ,\sigma _3W(\omega ^1) R _{\sigma _3(-\Delta +
  \omega ^1 )}( \overrightarrow{\mathbf{e}}\cdot  \alpha +\im \varepsilon)A _{\alpha 0}^{(1)}   \rangle \\& = \lim _{\varepsilon \to 0^+}\langle  \overline{A _{\alpha 0}^{(1)}} ,W(\omega ^1)^*\sigma _3W(\omega ^1) R _{\sigma _3(-\Delta +
  \omega ^1 )}( \overrightarrow{\mathbf{e}}\cdot  \alpha +\im \varepsilon)A _{\alpha 0}^{(1)}   \rangle .
\end{align*}
Using  the identity $W(\omega ^1)^*\sigma _3W(\omega ^1)=\sigma _3 Z(\omega ^1) W(\omega ^1)=\sigma _3 $, we conclude that
\begin{align*}
&     \langle \overline{A} _{\alpha 0} ,\sigma _3 R ^{+}_{ \mathcal{H}_{\omega ^1}  }( \overrightarrow{\mathbf{e }}\cdot  \alpha ) A _{\alpha 0}   \rangle   =  \langle  \overline{A _{\alpha 0}^{(1)}} , \sigma _3  R ^+_{\sigma _3(-\Delta +
  \omega ^1 )}(\overrightarrow{ \mathbf{e }}\cdot  \alpha  )A _{\alpha 0}^{(1)}   \rangle  \\& =  \langle  \overline{A _{\alpha 0}^{(1)}} , \sigma _3   P.V. \frac{1}{\sigma _3(-\Delta +
  \omega ^1 )  - \overrightarrow{\mathbf{e }}  \cdot  \alpha}    A _{\alpha 0}^{(1)}   \rangle   +\pi \im
  \langle  \overline{A _{\alpha 0}^{(1)}} , \sigma _3  \delta \({\sigma _3(-\Delta +
  \omega ^1 )  - \overrightarrow{\mathbf{e }}  \cdot  \alpha}\)      A _{\alpha 0}^{(1)}   \rangle   .
\end{align*}
In the last line, the first term is real valued and the last is imaginary. Hence, when we multiply \eqref{eq:last eq z}  by $ {\mathbf{e}}_j\overline{z}_j$, sum up on $j$ and take the imaginary part,   we obtain
\begin{equation*}
\begin{aligned}
    & 2^{-1}\partial _t \sum _{j }    {\mathbf{e}}_j
 | z _j|^2  =    -\pi  \sum_{\overrightarrow{\mathbf{e}}\cdot  \alpha     > \omega ^1} |z ^\alpha | ^2
     \langle  \overline{A _{\alpha 0}^{(1)}} , \sigma _3  \delta \({\sigma _3(-\Delta +
  \omega ^1 )  - \overrightarrow{\mathbf{e}}  \cdot  \alpha}\)      A _{\alpha 0 }^{(1)}   \rangle   \\& =  -\pi  \sum_{ \overrightarrow{\mathbf{e}} \cdot  \alpha     > \omega ^1} |z ^\alpha | ^2
   \left \langle  \begin{pmatrix} \overline{( A _{\alpha 0 }^{(1)})}_1  \\
\overline{( A _{\alpha 0}^{(1)})}_2
 \end{pmatrix}     ,  \begin{pmatrix} \delta ( -\Delta +\omega ^1   -  \overrightarrow{\mathbf{e}}    \cdot  \alpha )   &
0  \\
0 & -\underbrace{\delta (  \Delta -\omega ^1   -\overrightarrow{\mathbf{e}} \cdot  \alpha )}_{0}
 \end{pmatrix}        \begin{pmatrix} {(A _{\alpha  0}^{(1)})}_1  \\
{( A _{\alpha 0 }^{(1)})}_2
 \end{pmatrix}   \right  \rangle
 \\& = -\pi  \sum_{ \overrightarrow{\mathbf{e}} \cdot  \alpha     > \omega ^1} |z ^\alpha | ^2
   \left \langle \overline{( A _{\alpha 0 }^{(1)})}_1 , \delta ( -\Delta +\omega ^1   -  \overrightarrow{\mathbf{e}}   \cdot  \alpha ) {( A _{\alpha 0 }^{(1)})}_1 \right  \rangle ,
     \end{aligned}
\end{equation*}
where $^{t}{A _{\alpha 0 }^{(1)}} =\( {(A _{\alpha  0}^{(1)})}_1 , {(A _{\alpha 0 }^{(1)})}_2  \)$ and where
\begin{align}
&    \left \langle \overline{( A _{\alpha 0 }^{(1)})}_1 , \delta ( -\Delta +\omega ^1   -  \overrightarrow{\mathbf{e}}    \cdot  \alpha ) {( A _{\alpha 0 }^{(1)})}_1 \right  \rangle = \frac{1}{2(\overrightarrow{\mathbf{e}}    \cdot  \alpha - \omega ^1)}
\int _{|\xi |=\overrightarrow{\mathbf{e}}    \cdot  \alpha - \omega ^1 } \left | \widehat{(A _{\alpha  0}^{(1)})} _1 (\xi )  \right | ^2 dS\ge 0. \label{eq:FGR000}
\end{align}
The following is an hypothesis.
\begin{itemize}

 \item[(H8)]  We assume the last inequality to be strict, for appropriate choice of multi--indexes $\alpha$.

\end{itemize}
Then the argument closes up .

\begin{remark} \label{rem:check_FGR}   The coefficients in \eqref{eq:FGR000} are obtained after the NLS undergoes a significant number of coordinate
changes. As a consequence,   it is  not easy to write concretely   and   check numerically (H8). Notice though that in in \cite{CM20} the argument is much simplified, there is no normal forms argument and  the coefficients of the FGR are much simpler.
 \end{remark}

\section{ Further remarks and references  }
\label{sec:remarks}

We add some further remarks.

\begin{remark} \label{rem:soliton_well}    The problem of  the eventual behavior of a soliton of \eqref{NLS} in a confining well obtained adding to \eqref{NLS} a potential,
is mostly open.  For a non complete list of references see   \cite{BM16CMP,Bonanno15JDE,soffer, FGJS04CMP,HZ08IMRN,HZ1,Datchev,JLFGS06AHP},   and see therein for further references.
These papers treat long time behavior, but not asymptotic behavior. This problem is very similar to the oscillations discussed in Remark \ref{rem:numerical2}.
 \end{remark}

\begin{remark} \label{rem:non trapped soliton_well}     The effect of a potential on an escaping  soliton    is  easier to track, because, while it is deviated,
the soliton is preserved. There are various papers on the asymptotic behavior of  escaping potentials like \cite{CM1,CM15DCDS,deng,naum_raphael_escape_2018}.  A very suggestive
analysis of a soliton of the cubic integrable NLS in dimension 1 hitting a defocusing  delta potential is in Holmer, Marzuola and Zworski  \cite{HMZ1,HMZ2}. But the discussion,
 which uses also   the integrable structure, involves finite times only: it is not clear how to show that certain terms, that in \cite{HMZ2} are remainder, do not   develop in significant ones
 over larger intervals of time. Very little, beyond Deift and Zhou \cite{DZ},   is known about the use of the inverse scattering transform and the nonlinear steepest descent method    in the context of non--integrable systems,
 with the problem  in  \cite{HMZ2,MMS2} looking like    natural  for such a theory. For example, it would be natural to use the nonlinear steepest descent method to show that
 all solutions of a defocusing NLS in dimension 1 with a repulsive Dirac potential and with initial datum in $H^1(\R ) \cap \< x \> ^{-1}L^2(\R )$ decay like $t^{-1/2}$ and have the asymptotic profile
 that in  Masaki, Murphi and Segata \cite{MMS2} is proved only for small initial data.

 \noindent An asymptotic analysis over all times for a  problem  similar  to     \cite{HMZ2} is in Perelman  \cite{perelman09}, which however discusses a very flat soliton. In  \cite{perelman11} there is a finite time analysis of interaction of two solitons.

 \noindent Substantial modifications of a moving soliton in the presence of a NLS with a slowly varying coefficient in front of the nonlinearity are in  \cite{Munoz1,Munoz2}.

 \end{remark}

\begin{remark} \label{rem:blow up}      There are deep connections between the Fermi Golden Rule discussed here and the problem of the $\| \nabla u(t)\| _{L^2(\R ^d)}\sim \sqrt{\frac{\log |\log  t |}{t}}$  blow up
in the NLS with $\beta  (|u|^2) u=-|u|^{\frac{4}{d}}u$.  The key in the proof is the extension  in \cite{perelman01,MR4} of  the solitons in a larger class of functions, which are only approximately solutions of the NLS
and resemble the part of the solution here given by
\begin{equation*}
     u_{\text{approx}} = e^{\im \sigma _3  \tau \cdot \Diamond} \(  \Phi _{p } + P(p) \(  \sum_{l=1,..., {{n}}}z_l \xi_l   +
\sum_{l=1,..., {{n}}}\overline{z}_l  \sigma _1{{\xi }}_l  \)\) ,
\end{equation*}
obtained omitting the contribution of the continuous coordinate $f$. In \cite{MR4} there is the discussion of a Lyapunov function that takes the role of  $E _{\text{discr}} :=\sum _j \mathbf{e} _j |z_j|^2$.
In  \cite{MR4} the discussion is rather delicate because the coupling responsible for the Fermi Golden Rule is exponentially small, rather than polynomially small like in  Sect. \ref{subsec:FGR z}.
In \cite{perelman01,MR4} the choice of  \textit{profiles} (that is a generic solution is represented
as a sum of a profile plus a remainder, where the profile is similar to a ground state)  does not require an  algorithmic procedure (normal forms) and is related to the framework
introduced by previous authors and discussed by Sulem and Sulem \cite[Chapter 7]{Sulem}. Obviously, the step by step
normal forms argument glimpsed in Sect. \ref{sec:eigenvalues}, which eliminates   resonant terms one  at a time from the Taylor expansion of the energy $E$, would never yield the exponentially small resonant terms in  \cite{MR4}.  It is interesting that in the first paper of their series  \cite{MR1}, Merle and Raphael in Sect. 4.3 perform an argument similar to  normal forms which yields the non sharp upper bound on blow up
of \cite{MR1}. Presumably further changes of variables would yield algebraic improvements, which  nonetheless are not sharp. In any case,
from \cite{MR2} (the second paper in the series) on,  they settle in the optimal coordinates obtaining the sharp upper bound on blow up.
Notice that  in \cite{CM20}, in analogy to \cite{perelman01,MR4}, there is a choice of   profile that allows to avoid a normal forms argument.

Very delicate, especially because it is very difficult to estimate  in a sharp way various remainders,   is the proof of the   sharp lower bound on blow up in \cite{MR4},     where
a Lyapunov function is defined  starting from the  \textit{local virial identity} (the latter is stated in Proposition 2 \cite{MR4}, see also  the earlier \cite{MR1,MR2})
and then by various     adjustments.
The discussion is different from the one in the present survey, where the $E _{\text{discr}}$ is defined using the the discrete coordinates, which lose energy   leaking in the background.
Work needs to be done to compare the  Lyapunov function in \cite{MR4} with the $E _{\text{discr}}$ of  the present paper.
  It would  be interesting  to compare and unify the methods,   considering   problems mixing  the frameworks in  \cite{perelman01,MR4} and in here. One such problem might be  the one discussed in \cite{CM19JDE} and, by analogy, probably also problems   involving solitons trapped in wells, and in general, problems where the    linearizations have eigenvalues close to 0.
   Other similar problems are the ones involving the complicated      patterns in  \cite{MW10DCDS,GMW15DCDS}    near the  bifurcating standing waves of \cite{KirrKevPel,KKSW,kirrNatarajan}.

 \end{remark}

\begin{remark}
	Another topic which is not well studied is the relation between Fermi Golden Rule and the small ``wings" of nanopteron/micropterons \cite{Boyd98Book}.
	Here, a nanopteron/micropterons are infinite energy solutions which look  like   solitons locally but have a small nondecaying (or slow decaying) tail near spatial infinity. For mathematical results on the existence of such solutions for various equations see \cite{Beale91CPAM, HW17PD, JW20SAM, Lustri20PD, LP18SIAMJADS, Sun91JMAA}.
	When the tail is exponentially small w.r.t.\ a small parameter it is called nanopteron and, if it is polynomially small, it is called micropteron.
	Since the asymptotic stability result reviewed in this paper claims that there are no finite energy quasi-periodic solutions near solitons even though the linearized equation posses quasi-periodic solutions due to the internal modes, it is natural to ask if there exist   infinite energy quasi-periodic solutions near solitons, which should be   micropterons if they exist.
	Moreover, it is natural to guess that the ``wing" of such micropteron are related tothe  Fermi Golden Rule, in particular the first two terms of the r.h.s.\ of \eqref{eq:fermi2}.
	In   connection of Merle-Raphael's result on the blow up of critical NLS \cite{MR2,MR4}, it was shown by Johnson and Pan \cite{JP93PRSE} that there exists infinite energy solution which blows-up without the $\log\log$ correction.
	From the above point of view, it is natural to ask the relation between Merle and Raphael's optimal choice of the coordinate and Johnson-Pan's solution.
	Moreover the exponentially small Fermi Golden Rule and the asymptotic behavior of Johnson-Pan's solution at spatial infinity are of interest.
	However, these topics are completely open as far as the authors know.
	
\end{remark}

\begin{remark} \label{rem:small_sol}     Many of the papers on asymptotic stability of standing waves, focus on small standing waves which bifurcate from eigenvalues of a Schr\"odinger operator,
see \cite{SW1,SW2,PW}, \cite{TY1}-- \cite{TY4}, \cite{T,GNT,zhousigal,Gz,zhouweinstein1,zhouweinstein2,NPT,MMS1,CM19SIMA}. These papers impose restrictive hypotheses on the spectrum of the Schr\"odinger operator. Most of the hypotheses are dropped in \cite{CM15APDE}, which however treats only the case when the eigenvalues have multiplicity 1.  Higher multiplicities, but under restrictive conditions
 on the spectrum, are considered  in \cite{GP}.    Analogues of    \cite{CM15APDE} are for the NLKG   in    \cite {CMP2016}  and  for Dirac in  \cite{cuccagnatarullidirac}. Notice that in \cite{SW3,bambusicuccagna}, which treat
 NLKG, there are no standing waves because only real valued solutions are considered. In   \cite {CMP2016}, since complex valued solutions of the NLKG are considered, the dynamics
 of small energy solutions of the NLKG are more complicated than in \cite{SW3,bambusicuccagna}.
 \end{remark}

\begin{remark}
	The radiation damping also plays a role in the instability of excited states which are linearly stable.
	This   mechanism was called ``radiation induced instability" in \cite{HBW03SIMA} following  the name ''dissipation induced instability" \cite{BKMR94AIHPAN}.
	See also, \cite{CM16JNS}.
\end{remark}

\begin{remark} Global asymptotic results have been proved for equations where the nonlinearity is concentrated in a point, or in finitely many points, that is $\beta  (|u|^2) u $ is replaced by $ \delta (x-x_0)\beta  (|u|^2) u   $,
or by a linear combination of such terms. See \cite{Kom03}--\cite{KK10b},  \cite{Comech0} and therein.

Another model with very remarkable results is the energy critical focusing wave equation in 3 D, especially in the radial case, see \cite{Duyckaerts}, where the proof is based on the \textit{channel of energy} inequality,
which is specific to wave equations,
and on \textit{nonlinear profile} decompositions. In the context of the NLS, the   nonlinear profile  decompositions are rather complicated, see  \cite{NakanishiJMSJ}, and the presence, in the terminology of  \cite{PKA98}, of \textit{internal modes} of the solitons might render difficult proving the soliton decoupling, see also \cite{CM18IMRN}.
	
\end{remark}

\begin{remark}
	Little seems to be known about the nonlinear Klein Gordon Equations (NLKG). We do not know of any result analogous to  Theorems  \ref{thm:main rough},  \ref{thm:BP1} or \ref{thm:Cu1}  for
solitary waves of the NLKG.  Notice that an analogue of  Theorem   \ref{thm:main rough}  is known for   solutions with appropriate symmetries  of nonlinear Dirac Equations, see \cite{boussaidcuccagna}.
\end{remark}

\section*{Acknowledgments} S.C. was supported by a grant FRA 2018 from the University of Trieste.
M.M. was supported by the JSPS KAKENHI Grant Number 19K03579, JP17H02851 and JP17H02853.

Department of Mathematics and Geosciences,  University
of Trieste, via Valerio  12/1  Trieste, 34127  Italy.
{\it E-mail Address}: {\tt scuccagna@units.it}

Department of Mathematics and Informatics,
Faculty of Science,
Chiba University,
Chiba 263-8522, Japan.
{\it E-mail Address}: {\tt maeda@math.s.chiba-u.ac.jp}


\begin{thebibliography}{CP03}


\bibitem{adami13JMP}
R. Adami, D.Noja and C.Ortoleva, \emph{Orbital and asymptotic stability for standing waves of a nonlinear Schr\"odinger equation with concentrated nonlinearity in dimension three}, J. Math. Phys.  {54} (2013),
  no. 1, 013501, 33 pp.

\bibitem{asad}
R. Asad and G. Simpson,  \emph{ Embedded eigenvalues and the nonlinear Schr\"odinger equation,} J. Math. Phys.
52(3), 033511, 26 (2011)


\bibitem{Bambusi13CMPas}
D. Bambusi, \emph{Asymptotic stability of ground states in some
  {H}amiltonian {PDE}s with symmetry}, Comm. Math. Phys. \textbf{320} (2013),
  no.~2, 499--542.




\bibitem{bambusicuccagna}
  D.Bambusi and  S.Cuccagna, {\em On dispersion of
small energy solutions of the nonlinear Klein Gordon equation with a
potential},   Amer. Math. Jour., 133 (2011), 1421--1468 .




\bibitem{BM16CMP}
D.~Bambusi and A.~Maspero, \emph{Freezing of energy of a soliton in an external
  potential}, Comm. Math. Phys. \textbf{344} (2016), no.~1, 155--191.

\bibitem{Beale91CPAM}
J.~Thomas Beale, \emph{Exact solitary water waves with capillary ripples at
	infinity}, Comm. Pure Appl. Math. \textbf{44} (1991), no.~2, 211--257.

  \bibitem{beceanu}
  M. Beceanu,  {\em New estimates for a time-dependent Schr\"odinger equation}, Duke Math. Jour.   159  (2011), 417--477.

\bibitem{BKMR94AIHPAN}
A.~M. Bloch, P.~S. Krishnaprasad, J.~E. Marsden, and T.~S. Ratiu,
\emph{Dissipation induced instabilities}, Ann. Inst. H. Poincar\'{e} Anal.
Non Lin\'{e}aire \textbf{11} (1994), no.~1, 37--90.


\bibitem{Bonanno15JDE}
C. Bonanno, \emph{Long time dynamics of highly concentrated solitary waves
  for the nonlinear {S}chr\"{o}dinger equation}, J. Differential Equations
   {258} (2015),  717--735. 



\bibitem{borghese}
M. Borghese, R. Jenkins and   K. D. T.-R. McLaughlin,   \emph{
Long time asymptotic behavior of the focusing nonlinear Schr\"odinger equation},
Ann. Inst. H. Poincar\'e Anal. Non Lin\'eaire 35 (2018),   887--920.


\bibitem{boussaidcuccagna}N. Boussaid and S. Cuccagna,
\emph{On Stability of Standing Waves of Nonlinear Dirac Equations}, Communications in Partial Diff. Eq.   37 (2012), 1001--1056.



\bibitem{Boyd98Book}
J.~P. Boyd, \emph{Weakly nonlocal solitary waves and beyond-all-orders
	asymptotics}, Mathematics and its Applications, vol. 442, Kluwer Academic
Publishers, Dordrecht, 1998, Generalized solitons and hyperasymptotic
perturbation theory.

\bibitem{BG01MCS}
V.~S. Buslaev and V.~E. Grikurov, \emph{Simulation of instability of bright
  solitons for {NLS} with saturating nonlinearity}, Math. Comput. Simulation
  \textbf{56} (2001), no.~6, 539--546. 


\bibitem{Buslaev2008} V.Buslaev A. Komech, A.Kopylova  and D.Stuart,  {\em On asymptotic stability of solitary waves in Schr\"odinger quation coupled to nonlinear oscillaton},  Commun. Partial Differ. Equ., 33 (200),    669--705.


\bibitem{BP1}
V.Buslaev and  G.Perelman, {\em Scattering for the nonlinear
Schr\"odinger equation: states close to a soliton\/},  St.
Petersburg Math.J., 4
 (1993), 1111--1142.




\bibitem{BP2}
V.Buslaev and  G.Perelman, {\em On the stability of solitary waves for
nonlinear Schr\"odinger equations\/},   Nonlinear evolution
equations, editor N.N. Uraltseva, Transl. Ser. 2, 164, Amer. Math.
Soc.,
     75--98, Amer. Math. Soc., Providence (1995).


    \bibitem{BS}   V.Buslaev and  C.Sulem, {\em On the asymptotic
stability of solitary waves of Nonlinear Schr\"odinger equations},
Ann. Inst. H. Poincar\'e. An. Nonlin.,   {20} (2003),  419--475.


\bibitem{caz} {  T.~Cazenave},  {\em Semilinear Schr\"odinger equations,}Courant Lecture Notes in Mathematics, vol. {10}, Courant Lecture Notes, American Mathematical Society, Providence   (2003).


\bibitem{Chang}
S.M. Chang, S. Gustafson, K. Nakanish and T.P. Tsai,  {\em
Spectra of linearized operators for NLS solitary waves},
SIAM J. Math. Anal. 39 (2007/08),  1070--1111.

\bibitem{Comech0} A. Comech    {\em  Solutions with compact time spectrum to nonlinear Klein--Gordon and Schroedinger equations and the Titchmarsh theorem for partial convolution},
       Arnold Math. J. 5 (2019),   315--338.





\bibitem{Comech} A. Comech and  S. Cuccagna,  {\em On asymptotic stability of ground states of some systems of nonlinear Schr\"odinger equations},
       arXiv:1801.04079.


\bibitem{CP03CPAM}
A. Comech and D. Pelinovsky, \emph{Purely nonlinear instability of
  standing waves with minimal energy}, Comm. Pure Appl. Math. {56}
  (2003), 1565--1607. 






\bibitem{Costin}O.Costin, M.Huang and W. Schlag,  {\em
On the spectral properties of $L_{\pm}$ in three dimensions}
Nonlinearity 25 (2012),   125--164.







\bibitem{cuccagna2004}S. Cuccagna, \emph{A survey on asymptotic stability of ground states of nonlinear Schr\"odinger equations}, Dispersive nonlinear problems in mathematical physics, 21--57, Quad. Mat., 15, Dept. Math., Seconda Univ. Napoli, Caserta (2004).





\bibitem{Cuccagna11CMP}
S. Cuccagna, \emph{The {H}amiltonian structure of the nonlinear {S}chr\"odinger
  equation and the asymptotic stability of its ground states}, Comm. Math.
  Phys.  {305} (2011), no.~2, 279--331. 

\bibitem{Cuccagnatrans}
S. Cuccagna, \emph{On asymptotic stability of moving ground states of the nonlinear  {S}chr\"odinger
  equation}, Trans. Amer. Math. Soc.  {366} (2014), no.~6, 2827--2888.


\bibitem{Cu1}
  S.Cuccagna, {\em Stabilization of solutions to  nonlinear
Schr\"odinger equations},  Comm. Pure App. Math.  {54} (2001), pp.
1110--1145, erratum  Comm. Pure Appl. Math.  58  (2005), p. 147.



\bibitem{Cu2}   S.Cuccagna, {\em On asymptotic stability
of ground states of NLS},  Rev. Math. Phys.  {15} (2003),
877--903.





\bibitem{CM1} S.Cuccagna, M.Maeda, {\em On weak interaction between  a ground state and a non--trapping  potential},   J. Differential Eq.  {  256} (2014),   1395--1466.


\bibitem{CM15APDE}
S. Cuccagna and M. Maeda, \emph{On small energy stabilization in the
  {NLS} with a trapping potential}, Anal. PDE  {8} (2015),
  1289--1349. 


\bibitem{CM15DCDS}
S. Cuccagna and M. Maeda, \emph{On weak interaction between a ground state and a trapping potential}, Discrete Contin. Dyn. Syst. 35 (2015),
 3343--3376.


\bibitem{CM16JNS}
S. Cuccagna and M. Maeda, \emph{On orbital instability of spectrally
	stable vortices of the {NLS} in the plane}, J. Nonlinear Sci.  {26}
(2016), no.~6, 1851--1894.


\bibitem{CM18IMRN}
S. Cuccagna and M. Maeda, \emph{On Nonlinear Profile Decompositions and Scattering for an NLS--ODE Model}, International Mathematics Research Notices,  https://doi.org/10.1093/imrn/rny173.



\bibitem{CM19JDE}
S. Cuccagna and M. Maeda, \emph{Long time oscillation of solutions of nonlinear Schr\"odinger equations near minimal mass ground state}, Jour. of Diff. Equations, in press, https://doi.org/10.1016/j.jde.2019.11.047.


\bibitem{CM19SIMA}
S. Cuccagna and M. Maeda, \emph{On stability  of small solitons of the 1--D NLS with a  trapping delta potential}, SIAM J. Math. Anal.  51 (2019), 4311--4331.

\bibitem{CMP2016}
S. Cuccagna , M. Maeda Masaya and  Tuoc V. Phan,{ \em  On  small energy stabilization in the NLKG with a trapping potential}, Nonlinear Analysis,  146 (2016),  32--58.



\bibitem{CM20}
S. Cuccagna and M. Maeda, \emph{Coordinates at small energy and refined profiles for the Nonlinear Schr\" odinger Equation},  arXiv:2004.01366.



\bibitem{cuccagnamizumachi}
S.Cuccagna and  T.Mizumachi, {\em On asymptotic stability in energy
space of ground states for Nonlinear Schr\"odinger equations\/},
Comm. Math. Phys., 284
 (2008), pp.  51--87.

\bibitem{cupe2014}
S. Cuccagna and D. E. Pelinovsky,  {\em  The asymptotic stability of solitons in the cubic NLS equation on the line}, Appl. Anal. 93 (2014) 791--822.



\bibitem{CPV} S.Cuccagna, D.Pelinovsky and  V.Vougalter, {\em
Spectra of positive and negative energies in the linearization of
the NLS problem},  Comm.  Pure Appl. Math.  {58} (2005),  1--29.


\bibitem {cuccagnatarullidirac} S. Cuccagna and M. Tarulli, { \em On   stabilization   of small solutions   in the nonlinear Dirac equation with a trapping potential},   Journal of Mathematical Analysis and Applications,   436 (2016), 1332--1368.




\bibitem {cuccagnatarullidim2} S.Cuccagna and   M.Tarulli, {\em On asymptotic stability in energy space of ground states
of NLS in 2D} , Ann. I. H. Poincar\'e   26 (2009), 1361--1386.


\bibitem {cuccagnatarulliDNLS} S.Cuccagna and   M.Tarulli, {\em On asymptotic stability of standing waves of discrete
Schr\"odinger equation in $\Z$} , SIAM Jour. Math. An.   41 (2009), 861--885.



\bibitem{Datchev}
K. Datchev and J. Holmer,  {\em  Fast soliton scattering by attractive delta impurities,} Comm. Partial Differential Equations 34 (2009),   1074--1113.

\bibitem{BGR15Survey}
S. De~Bi\`evre, F. Genoud, and S. Rota~Nodari, {\em  Orbital
	stability: analysis meets geometry}, Nonlinear optical and atomic systems,
Lecture Notes in Math., vol. 2146, Springer, Cham, 2015,  147--273.

\bibitem{DZ} P. Deift and X.  Zhou,  \emph{  Perturbation theory for infinite-dimensional integrable systems on the line. A case study},
    Acta Math.
     188  (2002), 163--262.


\bibitem{Demanet}
  L. Demanet and W.  Schlag,    {\em
Numerical verification of a gap condition for a linearized nonlinear Schr\"odinger equation},
Nonlinearity 19 (2006),  829--852.


\bibitem{deng}
 Q. Deng, A. Soffer and  X. Yao, {\em Soliton-potential interactions for nonlinear Schr\"odinger equation in $\R^3 $}, SIAM J. Math. Anal. 50 (2018),   5243--5292.




\bibitem{Duyckaerts} T. Duyckaerts, C. Kenig  and F. Merle, {\em  Classification of radial solutions of the focusing, energy--critical wave equation}, Camb. J. Math. 1 (2013), 75--144.




\bibitem{soffer}  V. Fleurov and  A. Soffer, \emph{Soliton in a Well. Dynamics and Tunneling}, arXiv:1305.4279v1.



\bibitem{FGJS04CMP}
J.~Fr\"{o}hlich, S.~Gustafson, B.~L.~G. Jonsson  and I.~M. Sigal,
  \emph{Solitary wave dynamics in an external potential}, Comm. Math. Phys.
   {250} (2004),   613--642. 



\bibitem{galtabiar}
A.  Galtbayar and K. Yajima,  {\em The Lp-continuity of wave operators for one dimensional Schr\"odinger operators}, J. Math. Sci. Univ. Tokyo 7 (2000),  221--240.

\bibitem{Gz}  Zhou Gang, {\em Perturbation Expansion and N-th
Order Fermi Golden Rule of the Nonlinear Schr\"odinger Equations
\/},  J. Math. Phys., {48}( 2007), p. 053509.


\bibitem{zhousigal}
Zhou Gang and I.M.Sigal, {\em Relaxation of Solitons in Nonlinear
  Schr\"odinger Equations with Potential \/}, Advances in
  Math.,
  216 (2007),  443-490.


\bibitem{zhousigal2}Zhou Gang and I.M.Sigal, {\em Asymptotic stability of nonlinear
  Schr\"odinger Equations with Potential \/},  Reviews in Mathematical Physics, 17 (2005), 1193-1207.

\bibitem{zhouweinstein1}
Zhou Gang and M.I.Weinstein, {\em Dynamics of Nonlinear
Schr\"odinger/Gross-Pitaeskii Equations; Mass transfer in Systems
with Solitons and Degenerate Neutral Modes},  Anal. PDE  1 (2008),
pp.  267--322.

\bibitem{zhouweinstein2}
Zhou Gang and M.I.Weinstein, {\em Equipartition of Energy in Nonlinear
Schr\"odinger/Gross-Pitaeskii Equations},   Appl. Math. Res. Express. AMRX (2011),  123--181.





\bibitem{GMW15DCDS}
R.~H. Goodman, J.~L. Marzuola  and M.~I. Weinstein,
  \emph{Self-trapping and {J}osephson tunneling solutions to the nonlinear
  {S}chr\"odinger/{G}ross-{P}itaevskii equation}, Discrete Contin. Dyn. Syst.
  {35} (2015), 225--246. 





\bibitem{G1} M.Grillakis,  {\em  Analysis of the linearization
around a critical point of an infinite dimensional Hamiltonian
system},  Comm. Pure Appl. Math    43 (1990), 299--333



\bibitem{GNT}
S.Gustafson,   K.Nakanishi and  T.P.Tsai,  {\em Asymptotic stability and completeness in the energy space for nonlinear Schr\"odinger equations with small solitary waves}, Int. Math. Res. Not.  {   66} (2004)  , 3559--3584.



\bibitem{GP}
S.Gustafson,    T.V.Phan,  {\em Stable directions for degenerate excited states of nonlinear  Schr\"odinger equations}, SIAM J. Math. Anal. {  43} (2011) , 1716--1758.


\bibitem{HW17PD}
A. Hoffman and J.~Douglas Wright, \emph{Nanopteron solutions of diatomic
	{F}ermi-{P}asta-{U}lam-{T}singou lattices with small mass-ratio}, Phys. D
\textbf{358} (2017), 33--59.

\bibitem{HBW03SIMA}
P. Hagerty, A.~M. Bloch, and M.~I. Weinstein, \emph{Radiation
	induced instability}, SIAM J. Appl. Math.  {64} (2003/04),
484--524.

\bibitem{HMZ1}
J. Holmer, J Marzuola  and M. Zworski, \emph{  Soliton splitting by external delta potentials}, J. Nonlinear Sci. 17 (2007), no. 4, 349--367.


 \bibitem{HMZ2} J. Holmer, J Marzuola  and M. Zworski, \emph{  Fast soliton scattering by delta impurities}, Comm. Math. Phys. 274 (2007),   187--216.






\bibitem{HZ08IMRN}
J. Holmer and M. Zworski, \emph{Soliton interaction with slowly varying
  potentials}, Int. Math. Res. Not. IMRN (2008),  Art. ID rnn026, 36.



	
		\bibitem {HZ1}
		J.Holmer,   M. Zworski,  {\em Slow soliton interaction with delta impurities}
J. Mod. Dyn. 1 (2007),   689--718.







\bibitem{how1}J. S.  Howland, {\em On the Weinstein-Aronszajn formula }, Arch. Rational Mech. Anal. 39 (1970),
323--339.


\bibitem{how2}J. S.  Howland, {\em Puiseux series for resonances at an embedded eigenvalue},  Pacific J. Math. 55
(1974), 157--176.

\bibitem{JW20SAM}
M.~A. Johnson and J.~Douglas Wright, \emph{Generalized solitary waves in
		the gravity--capillary   Whitham equation}, Studies in Applied Mathematics
\textbf{144} (2020), no.~1, 102--130.

\bibitem{JP93PRSE}
R. Johnson and X.~Bin Pan, \emph{On an elliptic equation related to the
	blow-up phenomenon in the nonlinear {S}chr\"{o}dinger equation}, Proc. Roy.
Soc. Edinburgh Sect. A \textbf{123} (1993), no.~4, 763--782.

\bibitem{JLFGS06AHP}
B.~L.~G. Jonsson, J. Fr\"{o}hlich, S. Gustafson, and
  I.~M. Sigal, \emph{Long time motion of {NLS} solitary waves in a
  confining potential}, Ann. Henri Poincar\'{e}  {7} (2006),
  621--660. 


\bibitem{JSS}
 J.L.Journ\' e, A.Soffer and  C.D.Sogge,
{\em Decay estimates for Schrodinger operators }, Comm.P. Appl. Mat.
  44 (1991),  573--604.




\bibitem{kato}
T.Kato, {\em Wave operators and similarity for some non-selfadjoint
operators \/},  Math. Annalen, 162
 (1966),
  258--269.


\bibitem{Kl-Tao} { M.~Keel and T,~Tao},
{\em Endpoint Strichartz estimates},  Amer. J. Math. {120} (1998),
955--980.


\bibitem{KirrKevPel} E.Kirr,  P.G.Kevrekidis and D.E. Pelinovsky.
{\em Symmetry breaking bifurcation in Nonlinear
Schr\"odinger equation with symmetric potentials},  Comm. Math. Phys., 308 (2011),  795--844.



\bibitem{KKSW} E.Kirr,  P.G.Kevrekidis, E.Shlizerman and  M.I.Weinstein.
{\em Symmetry breaking bifurcation in Nonlinear
Schr\"odinger/Gross-Pitaevskii Equations,}  SIAM J. Math. Anal. 40
(2008),. 566--604.



\bibitem{kirr1}
E. Kirr and O. Mizrak,  {\em Asymptotic stability of ground states in 3D nonlinear Schr \"odinger equation including subcritical cases}, J. Funct. Anal.257 (2009),  3691--3747



 \bibitem{kirrNatarajan}
E. Kirr and V. Natarajan,  {\em
    The global bifurcation picture for ground states in nonlinear Schr \"odinger equations},  arXiv:1811.05716.




 \bibitem{kirr2}
 E. Kirr and A. Zarnescu, {\em Asymptotic stability of ground states in 2D nonlinear Schr \"odinger equation including subcritical cases} , J. Differential Equations 247 (2009), 710--735.


\bibitem{Kom03}
  A.  Komech,  {\em On attractor of a singular nonlinear U(1)-invariant Klein-Gordon equation}, in Progress in analysis, Vol. I, II (Berlin, 2001),
pp. 599--611, World Sci. Publ., River Edge, NJ, 2003.


\bibitem{KK09}   A. Komech and A. Komech, {\em   Global attraction to solitary waves for Klein--Gordon equation with mean field interaction}, Ann. Inst. H.
Poincar\'e Anal. Non Lin\'eaire 26 (2009),   855--868.



\bibitem{KK10a}   A. Komech and A. Komech, {\em Global attraction to solitary waves for a nonlinear Dirac equation with mean field interaction}, SIAM J. Math.
Anal. 42 (2010),   2944--2964.


\bibitem{KK10b}   A. Komech and A. Komech, {\em  On global attraction to solitary waves for the Klein-Gordon field coupled to several nonlinear oscillators}, J.
Math. Pures Appl.   93 (2010),   91--111.

\bibitem{Komech2012}
A.Komech, E.Kopylova and  D.Stuart, {\em On asymptotic stability of solitons in a nonlinear Schr\"odinger equation\/}, Commun. Pure Appl. Anal. 11 (2012),   1063--1079.



\bibitem{KS} J.Krieger and  W.Schlag, {\em Stable manifolds for all
monic supercritical focusing nonlinear Schr\"odinger equations in
one dimension\/}, J. Amer. Math. Soc., 19 (2006), pp. 815--920.




\bibitem {KMM2}
M. Kowalczyk, Y. Martel and C. Mu\~{n}oz,
{\em  Kink dynamics in the $\phi ^4$ model: asymptotic stability for odd perturbations in the energy space},
J. Amer. Math. Soc. 30 (2017), 769--798.


\bibitem {KMM3}
M. Kowalczyk, Y. Martel and C. Mu\~{n}oz,
{\em  Soliton dynamics for the 1D NLKG equation with symmetry and in the absence of internal modes},  Jour. of the Europ. Math. Soc., to appear.


\bibitem {KMM4}M. Kowalczyk, Y. Marte,  C. Mu\~{n}oz and H. Van Den Bosch {\em
A sufficient condition for asymptotic stability of kinks in general $(1+1)$-scalar field models}, arXiv:2008.01276.


\bibitem{Lustri20PD}
C.~J. Lustri, \emph{Nanoptera and {S}tokes curves in the 2-periodic
	{F}ermi--{P}asta--{U}lam--{T}singou equation}, Phys. D \textbf{402} (2020),
132239.

\bibitem{LP18SIAMJADS}
C.~J. Lustri and M.~A. Porter, \emph{Nanoptera in a period-2 {T}oda
	chain}, SIAM J. Appl. Dyn. Syst. \textbf{17} (2018), no.~2, 1182--1212.

\bibitem {MM1}Y. Martel and F. Merle, {\em Asymptotic stability of solitons of the gKdV equations with general nonlinearity}, Math. Ann. 341 (2008),  391--427.







\bibitem {Marzuola}
J. Marzuola and G. Simpson, Gideon, {\em
Spectral analysis for matrix Hamiltonian operators}
Nonlinearity 24 (2011),   389--429.



\bibitem{MRS10JNS}
J.~L. Marzuola, S. Raynor, and G. Simpson, \emph{A system of {ODE}s
  for a perturbation of a minimal mass soliton}, J. Nonlinear Sci. {20}
  (2010),  425--461. 


\bibitem{MW10DCDS}
J.~L. Marzuola and M.~I. Weinstein, \emph{Long time dynamics near the
  symmetry breaking bifurcation for nonlinear
  {S}chr\"odinger/{G}ross-{P}itaevskii equations}, Discrete Contin. Dyn. Syst.
  {28} (2010),   1505--1554. 




\bibitem{MMS1} S. Masaki, J. Murphy and  J. Segata, \emph{Stability of small solitary waves for the 1$d$ NLS with an attractive delta potential}, 	to appear in Anal.  PDE..

\bibitem{MMS2}
S. Masaki, J. Murphy, and J. Segata,\emph{ Modified scattering for the one-dimensional cubic NLS
with a repulsive delta potential}, Int. Math. Res. Not. doi:10.1093/imrn/rny011 .





\bibitem{MR4}F.Merle and P.Raphael, \emph{On a sharp lower bound on the blow-up rate for the $L^2$ critical nonlinear Schr\"odinger equation }, J. Amer. Math. Soc.
   {19} (2006),  37--90.



\bibitem{MR2}
F.Merle and P.Raphael, \emph{Sharp upper bound on the blow--up rate for the critical nonlinear  {S}chr\"odinger equation}, Geometric \& Functional Analysis   {13} (2003),
    591--642.


\bibitem{MR1}
F.Merle and P.Raphael, \emph{The blow-up dynamic and upper bound on the blow-up rate for critical nonlinear  {S}chr\"odinger equation}, Ann. of Math.  {161} (2005),
    157--222.



\bibitem{M1}  T.Mizumachi, {\em Asymptotic stability of small
solitons to 1D NLS with potential }, Jour. of Math.   Kyoto
University, 48 (2008),  471-497.
%

\bibitem{M2}  T.Mizumachi, {\em Asymptotic stability of small solitons
for 2D Nonlinear Schr\"{o}dinger equations with potential}, Jour. of
Math.   Kyoto University,   43 (2007),  599-620.




\bibitem  {Munoz1}  C. Munoz, {\em Sharp inelastic character of slowly varying NLS solitons},   arXiv:1202.5807.


\bibitem  {Munoz2}  C. Munoz, {\em On the soliton dynamics under slowly varying medium for Nonlinear Schr\"odinger equations},     Math. Ann., 353 (2012), 867--943.



  \bibitem {NakanishiJMSJ}
  K.Nakanishi,  {\em Global dynamics below excited solitons for the nonlinear Schr\"odinger equation with a potential},    J. Math. Soc. Japan,  {69} (2017),   1353--1401.




\bibitem {NPT}K.Nakanishi, T.V.Phan and  T.P.Tsai, {\em Small solutions of nonlinear Schr\"odinger equations near first excited states}, Jour. Funct. Analysis  {  263}  (2012), 703--781.




\bibitem  {naum_raphael_escape_2018}
I. Naumkin, P. Raphael, {\em On travelling waves of the non linear  Schr\"odinger equation escaping a potential well},  Ann. Henri Poincar\' e 21 (2020),  1677--1758.

\bibitem{PKA98}
D.~E. Pelinovsky, Y.~S. Kivshar, and V.~V. Afanasjev,
\emph{Internal modes of envelope solitons}, Phys. D \textbf{116} (1998),
no.~1-2, 121--142.

\bibitem  {perelman11} G.Perelman, {\em  Two soliton collision for nonlinear Schr\"odinger equations in dimension 1},
Ann. Inst. H. Poincar\'e Anal. Non Lin\'eaire 28 (2011),   357--384.


 \bibitem  {perelman09} G.Perelman, {\em   A remark on soliton-potential interactions for nonlinear Schr\"odinger  equations},   Math. Res. Lett. 16 (2009),  477--486.

\bibitem  {perelman3} G.Perelman, {\em  Asymptotic stability of multi-soliton solutions for nonlinear     Schr\"odinger equations},  Comm. Partial Diff.  29 (2004),  1051--1095.


\bibitem{perelman01}
 G. Perelman, {\em On  the Formation of Singularities in Solutions of the Critical
Nonlinear Schr\"odinger equation},   Ann. Henri Poincar\'e 2 (2001), 605--673.






\bibitem{PW}
C. Pillet and E. Wayne, {\em Invariant manifolds for a class of
dispersive, Hamiltonian partial differential equations},  J. Diff. Eq.
  141 (1997),  310--326.

\bibitem  {RSS2} I.Rodnianski, W.Schlag and  A.Soffer, {\em
Asymptotic stability of N-soliton states of NLS}, preprint (2003),
	arXiv:math/0309114v1 .


\bibitem  {saalmann}
 A. Saalmann,   {\em Asymptotic stability of N-solitons in the cubic NLS equation},  J. Hyperbolic Differ. Equ. 14 (2017),  455--485.




\bibitem{Schlag}   W.Schlag,  {\em Stable manifolds for an orbitally
unstable NLS \/},  Ann. of Math.    169  (2009),  139--227.

\bibitem{Sigal93CMP}
I.~M. Sigal, \emph{Nonlinear wave and {S}chr\"odinger equations. {I}.
	{I}nstability of periodic and quasiperiodic solutions}, Comm. Math. Phys.
\textbf{153} (1993), no.~2, 297--320.

\bibitem{SW1} A.Soffer and  M.I.Weinstein, {\em  Multichannel nonlinear
scattering for nonintegrable equations \/}, Comm. Math. Phys., 133
(1990),   116--146.






\bibitem{SW2}
A.Soffer and  M.I.Weinstein, {\em  Multichannel nonlinear scattering II.
The case of anisotropic potentials and data \/},  J. Diff. Eq., 98
 (1992),
  376--390.


\bibitem{SW3}
A.Soffer and  M.I.Weinstein, {\em Resonances, radiation damping and
instability in Hamiltonian nonlinear wave equations \/},  Invent.
Math., 136
 (1999),
  9--74.






\bibitem{SW4}    A.Soffer and  M.I.Weinstein,
{\em Selection of the ground state for nonlinear Schr\"odinger
equations}, Rev. Math. Phys.  16 (2004),    977--1071.

\bibitem{Sulem}C.Sulem and P.-L. Sulem,  {\em The Nonlinear Schr\"odinger Equation}, Applied Mathematica Sciences vol. 139 (1999), Springer, New York.

\bibitem{Sun91JMAA}
S.~M. Sun, \emph{Existence of a generalized solitary wave solution for water
	with positive {B}ond number less than {$1/3$}}, J. Math. Anal. Appl.
\textbf{156} (1991), no.~2, 471--504.

\bibitem{Stuart08}
C.~A. Stuart, {\em Lectures on the orbital stability of standing waves and	application to the nonlinear {S}chr\"odinger equation}, Milan J. Math.
 {76} (2008), 329--399.


\bibitem{T}   T.P.Tsai, {\em  Asymptotic dynamics of nonlinear
Schr\"odinger equations with many bound states},    J. Diff. Eq.
 {192}  (2003),     225--282.


\bibitem{TY1}
  T.P.Tsai and  H.T.Yau, {\em Asymptotic dynamics of nonlinear
Schr\"odinger equations: resonance dominated and radiation dominated
solutions}, Comm. Pure Appl. Math.  {55}  (2002),   153--216.

\bibitem{TY2}
  T.P.Tsai and  H.T.Yau, {\em Relaxation of excited states in
nonlinear Schr\"odinger equations}, Int. Math. Res. Not.  {31}
(2002),   1629--1673.

\bibitem{TY3}
{  T.P.Tsai and  H.T.Yau}, {\em Classification of asymptotic profiles
for nonlinear Schr\"odinger equations with small initial data}, Adv.
Theor. Math. Phys.  {6} (2002),  107--139.

\bibitem{TY4}
{  T.P.Tsai and  H.T.Yau}, {\em   Stable directions for
excited states of nonlinear Schr\"odinger equations},  Comm.
Partial Diff. Eq.   27  ( 2002)  2363--2402.


\bibitem{weder1}R. Weder, {\em The $W^{k,p}$--Continuity of the Schr\"odinger wave operators on the line}, Comm.  Math. Physics  208  (1999), 507--520.


\bibitem{weder2}R. Weder, {\em $L^p\to L^{p^\prime}$ estimates for
 the Schr\"odinger equation
   on the line and inverse
scattering for the nonlinear Schr\"odinger equation with a potential},
   J. Funct. Anal.   170 (2000) , 37--68.





\bibitem{W2}
  M.I.Weinstein, {\em Modulation stability of ground states of
nonlinear Schr\"odinger equations},  Siam J. Math. Anal.  16 (1985),
 472--491.



\bibitem{Y1}   K.Yajima, {\em  The $W^{k,p}$-continuity of wave operators
for Schr\"{o}dinger operators}, J. Math. Soc. Japan,  {47} (1995),
   551--581.
%
\bibitem{Y2} K.Yajima, {\em  The $W^{k,p}$-continuity of wave operators
for Schr\"{o}dinger operators III.}, J. Math. Sci. Univ. Tokyo,
 {2} (1995),   311--346.



\end{thebibliography}
\end{document}